\documentclass{birkmult}
\usepackage{amsmath, amsfonts,amsthm, amssymb,graphicx}
\usepackage[usenames]{color}

\usepackage{amsmath, amsfonts,amsthm, amssymb,graphicx}
\usepackage[usenames]{color}

 \newtheorem{thm}{Theorem}[section]
 \newtheorem{cor}[thm]{Corollary}

 \newtheorem{prop}[thm]{Proposition}

 \newtheorem{ex}[thm]{Example}
 \newtheorem{defn}[thm]{Definition}
 \newtheorem{rem}[thm]{Remark}
 
 \newtheorem{assu}[thm]{Assumption}

 \newcommand{\zz}{\mathfrak{Z}}

\long\def\symbolfootnote[#1]#2{\begingroup%
\def\thefootnote{\fnsymbol{footnote}}\footnote[#1]{#2}\endgroup}

\begin{document}

\title[Approximate radical for clusters]{Approximate radical for clusters: a global approach using Gaussian elimination or SVD}

\author[I. Janovitz-Freireich]{Itnuit
Janovitz-Freireich}\address{Mathematics Department,\;North
Carolina State University, \;Raleigh, NC,
USA.}\email{ijanovi2@ncsu.edu}

\author[L. R\'{o}nyai]{Lajos R\'{o}nyai}\address{Computer and Automation Institute,\;
Hungarian Academy of Sciences and Budapest University of
Technology and Economics, Budapest,
Hungary.}\email{lajos@csillag.ilab.sztaki.hu}

\author[\'{A}. Sz\'{a}nt\'{o}] {\'{A}gnes
Sz\'{a}nt\'{o}}\address{Mathematics Department,\;North Carolina
State University, \;Raleigh, NC, USA.}\email{aszanto@ncsu.edu}
\date{\today}

\begin{abstract}
We present a method based on Dickson's lemma to compute the
``approximate radical" of a zero dimensional ideal $\tilde I$ in
$\mathbb{C}[x_1, \ldots, x_m]$ which has zero clusters: the
approximate radical ideal has exactly one root in each cluster for
sufficiently small clusters. Our method  is ``global" in the sense
that it does not require any local approximation of the zero
clusters: it reduces the problem to the computation of the
numerical nullspace of the so called ``matrix of traces", a matrix
computable from the generating polynomials of $\tilde I$. To
compute the numerical nullspace of the matrix of traces we propose
to use Gauss elimination with pivoting or singular value
decomposition. We prove that if $\tilde I$ has $k$ distinct zero
clusters each of radius at most $\varepsilon$ in the
$\infty$-norm, then $k$ steps of Gauss elimination on the matrix
of traces yields a submatrix with all entries
\text{asymptotically} equal to $\varepsilon^2$. We also show that
the $(k+1)$-th singular value of the matrix of traces is
proportional to $\varepsilon^2$. The resulting approximate radical
has one root in each cluster with coordinates which are the
arithmetic mean of the cluster, up to an error term asymptotically
equal to $\varepsilon^2$. In the univariate case our method gives
an alternative to known approximate square-free factorization
algorithms which is simpler and its accuracy is better understood.
\end{abstract}

\thanks{This work was completed with the support of NSF grants
CCR-0306406 and CCR-0347506 and OTKA grants T42481 and T42706 and
NK63066.}


\maketitle


\section*{Introduction}

Let $I\subset \mathbb{C}[{\bf x}]$ be a polynomial ideal in $m$
variables  ${\bf x}=[x_1,\ldots,x_m]$ with roots ${\bf z}_1, \ldots,
{\bf z}_k\in \mathbb{C}^m$ of multiplicities $n_1, \ldots, n_k$,
respectively, and let $\tilde{I}\in \mathbb{C}[{\bf x}]$ be an ideal
with clusters $C_1, \ldots, C_k$ such that each cluster $C_i$ has
$n_i$ roots around ${\bf z}_i$ within radius $\varepsilon$ in the
$\infty$-norm for $i=1, \ldots, k$. We present an algorithm which
computes an \emph{approximate radical} of $\tilde{I}$, denoted by
$\widetilde{\sqrt I}$, which has exactly one root for each cluster,
and we show that such root corresponds to the arithmetic mean of the
cluster.

The method we present in the paper is ``global" in the sense that
we do not  use any local information about the  roots in the
clusters, only  the coefficients of the system of polynomials
defining $\tilde{I}$, and we return another system of polynomials
where all near multiplicities are eliminated.  In the univariate
case such global algorithms are used for example in approximate
factoring (see \cite{KalMay2003}), where the input polynomial
needs to be ``square-free"  in the approximate sense. Previous
global methods which handle univariate polynomials with clusters
use approximate gcd computation and approximate polynomial
division in order to either  factor out the near multiplicities or
to compute the approximate multiplicity structure and find the
roots of the nearest polynomial with the given multiplicity
structure \cite{SasNod89,hriste97,KalMay2003,Zeng2003}. The method
we propose here offers an alternative algorithm to factor out near
multiplicities, which is simpler, and the relation between the
accuracy of the output and the size of the clusters is better
understood. We describe separately our method applied to the
univariate case, and illustrate its simplicity and accuracy.

Our method is based on Dickson's lemma, which gives the Jacobson
radical of a   finite dimensional associative algebra over a field
of characteristic $0$ via the vanishing of traces of elements in the
algebra. An immediate application of Dickson's lemma to the algebra
$\mathbb{C}[\bold{x}]/I$ finds a basis for $\sqrt{I}/ I$ by finding
the nullspace of the \emph{matrix of traces} $R$, a matrix
computable from the generating polynomials of $I$ using either
multiplication matrices or other trace computation methods, as
described below.

The main focus of the paper is to adapt the method based on
Dickson's lemma to the case when the ideal $\tilde{I}$ has
clusters of roots. In the paper we assume that  both
$\mathbb{C}[\bold{x}]/I$ and $\mathbb{C}[\bold{x}]/\tilde{I}$ are
finite dimensional over $\mathbb{C}$ and have the same basis
$B\subset\mathbb{C}[{\bf x}]$. Note that if $I$ is generated by a
well-constrained system, then ``almost all" perturbations
$\tilde{I}$ of $I$ will satisfy our assumption, however our
results are not limited to well-constrained systems only. On the
other hand, the results we prove in this paper measure the
accuracy of the output in terms of the size of the clusters, as
opposed to the size of the perturbation of the generating
polynomials of the ideal $I$. The extension of our method to
handle perturbations which change the structure of the factor
algebra and to understand the accuracy of the output in terms of
the size of the coefficient perturbation is the subject of future
research. The results in this paper  can be summarized as follows:

Given the basis $B$  and the matrix of traces $R$ associated to
$\tilde I$ and $B$,  using Gaussian elimination with complete
pivoting (GECP) we give asymptotic estimates of order
$\varepsilon^2$ for the ``almost vanishing" entries in $U_k$, the
partially row reduced matrix of $R$, as well as upper bounds for the
coefficients of $\varepsilon^2$, where $\varepsilon$ is the radius
of the clusters in the $\infty$-norm. These bounds can be used to
give a threshold to decide on the  numerical rank of $R$, and   to
indicate the relationship between the numerical rank and the size of
the clusters.

Alternatively, we show how our results for the GECP of the matrix
of traces $R$ imply asymptotic bounds on the singular values of
$R$. We also obtain in this case that the ``almost vanishing"
singular values are proportional to the square of the size of the
clusters. This implies that for the numerical rank determination
of $R$, computing its $SVD$ works similarly as using GECP.

 Using a basis of the numerical
nullspace of $R$ (or possibly an extended version of it), we
define a set of generating polynomials for the approximate radical
ideal $\widetilde{\sqrt I}$, or similarly, define a system of
multiplication matrices $M'_{x_1},\ldots,M'_{x_m}$  of
$\mathbb{C}[\bold{x}]/\widetilde{\sqrt I}$ with respect to a basis
$B'$. We prove that modulo $\varepsilon^2$ the generating
polynomials of $\widetilde{\sqrt I}$ are consistent and have roots
with coordinates which are the arithmetic means of the coordinates
of the roots in the clusters, which implies that  the matrices
$M'_{x_1},\ldots,M'_{x_m}$ commute and their eigenvalues are the
arithmetic means of the coordinates of the roots in the clusters,
all modulo $\varepsilon^2$. In other words, our algorithm finds
the coefficients of a polynomial system with  roots which are the
means of the clusters up to a precision of about twice as many
digits as the radius of the clusters, assuming that the clusters
are sufficiently small.

 Let us briefly mention  some of the possible methods to compute the matrix of traces
 $R$, although in the paper we do
not elaborate on this aspect.
 As we shall demonstrate in the paper, the matrix of traces $R$ is readily
computable from a system of multiplication matrices of
$\mathbb{C}[\bold{x}]/\tilde I$, for example from
$M_{x_1},\ldots,M_{x_m}$, where  $M_{x_i}$ denotes the
  matrix of the multiplication
map by  $x_i$ in $\mathbb{C}[\bold{x}]/\tilde I$ written  in terms
of the basis $B$. One can compute $M_{x_i}$ using Gr\"obner bases
(see for example \cite{Cor96}), resultant and subresultant
matrices \cite{mandem95,Cha,Sza2}, Lazard's algorithm
\cite{Laz,CoGiTrWa95}, or by methods that combine these
\cite{Mou05}. Thus, our algorithm   reduces the problem of finding
the eigenvalues of matrices $M_{x_1},\ldots,M_{x_m}$ which have
clustered eigenvalues   to finding eigenvalues of the smaller
matrices $M'_{x_1},\ldots,M'_{x_m}$ with well separated
eigenvalues.

 In certain cases,  the matrix of traces can be computed directly
from the generating polynomials of $\tilde I$, without using
multiplication matrices.  We refer to the papers \cite{brigonz02,
DiazGonz01,CatDicStu96,CatDicStu98, carmou96} for the computation of
traces using residues and Newton sums, or \cite{DaJe05} using
resultants.

 Also, fast computation techniques like the ``baby
steps-giant steps" method \cite{Ka92:issac, Shou99, Scho05} can be
implemented to speed up the computation of all $n^2$ entries of the
\text{matrix} of traces. As we  prove in the paper, the entries of
the matrix of traces $R$ are continuous in the size $\varepsilon$ of
the root perturbation around $\varepsilon=0$, unlike the entries of
multiplication matrices which may have many accumulation points as
$\varepsilon$ approaches zero. Therefore, avoiding the computation
of the multiplication matrices has the advantage of staying away
from the possible large computational errors caused by the
discontinuity of their entries.

In the multivariate case, most of the methods handling clusters of
roots in the literature are
 ``local" in that they
assume  sufficiently close approximations for the clusters in
question. Our algorithm, viewed as having the multiplication
matrices as input,  is closest to the approach in
 \cite{mandem95, corlessgiani97} in that these papers also
reduce the problem to  the computation of the eigenvalues of a
system of approximate multiplication matrices. Both of these
papers propose to reorder the eigenvalues  of the multiplication
matrices to group the clusters together. For the reordering of the
eigenvalues these papers compute approximations of the eigenvalues
by either using the approach  in \cite{BaDeMc89}
 or using the univariate method of \cite{hriste97}.
 In contrast, our method reorders the
eigenvalues of all multiplication matrices simultaneously without
approximating the eigenvalues, grouping one eigenvalue from each of
the clusters together in a way which facilitates the computation of
the means of the clusters and the elimination of the rest of the
nearly repeated eigenvalues. Another local method to handle near
multiple roots is the ``deflation'' algorithm, studied in the works
 \cite{OjWa1983,Oj1987a,Oj1987b,Lecerf02,LeVeZh05},   to replace the original system
which had a near multiple root with another one which has the same
root with multiplicity one, using an approximation of the root in
question. Related to the deflation algorithm, in
\cite{stetter96,Stetter04, DaZeng05} methods are proposed to compute
the multiplicity structure of a root locally in terms of the so
called dual basis, and then computing good approximations for the
individual roots in the cluster, assuming that either a near system
with multiple roots  is known, or a sufficient approximation of the
multiple root is given. Additionally, methods for computing singular
solutions of both polynomials and analytic functions using homotopy
continuation can be found in \cite{MoSoWa91, MoSoWa92, MoSoWa92b}.

We also include here reference to some of the related methods for
solving systems of polynomial equations with exact multiplicities:
involving the computation of dual bases \cite{moste95, mamomo95,
Stetter04}, or in the univariate (or bivariate) case,  using Gauss
maps \cite{kosakpa04}, or analyzing the structure of the
multiplication matrices by transforming them to an upper triangular
form \cite{yokonorotake92, moritsukuri97, moritsukuri99}. Previous
work using Dickson's Lemma to compute radical ideals in the exact
case includes \cite{ArSo95, becwor96}. Also, \cite{Rouiller99} uses
trace matrices in order to find separating linear forms
deterministically.

The present paper is the extended and unabridged   version of the
paper that appeared in \cite{JaRoSza06}.

\subsection*{Acknowledgements:} We would like to thank Erich Kaltofen
for suggesting the problem.

\section{Preliminaries}\label{pre}

Let $A$ be an associative algebra over a field $F$ of characteristic
$0$. (See definition and basic properties of associative algebras in
\cite{ronyai85,Pie}.)

An element $x\in A$ is \emph{nilpotent} if $x^m=0$ for some positive
integer $m$.

An element $x\in A$ is \emph{properly nilpotent} if $xy$ is
nilpotent for every $y\in A$.

The \emph{radical} of $A$, denoted $Rad(A)$, is the set of properly
nilpotent elements of $A$. The radical $Rad(A)$ is an ideal of $A$.
In commutative algebras  nilpotent elements are properly nilpotent,
hence for a commutative $A$ the radical $Rad(A)$ is simply the set
of nilpotent elements in $A$.

Throughout the paper we assume that $A$ is finite dimensional over
$F$. Fix a basis $B=[b_1,\ldots,b_n]$ of $A$ (note that later we
will need to fix the order of the elements in $B$, that is why we
use vector notation). We call the \emph{multiplication matrix}
$M_x$ of $x\in A$ the matrix of the multiplication map
\[
\begin{aligned}
m_x:\;A&\longrightarrow A\notag\\
      [g]&\;\mapsto [xg]\notag\\
\end{aligned}
\]
written in the basis $B$. It is easy to verify  (cf. Page 8 in
\cite{Pie}) that the map $x\mapsto M_x$ is an algebra homomorphism,
called {\em regular representation} from $A$ to $M_n(F)$.

The \emph{trace of x}, denoted $Tr(x)$, is the trace of the matrix
$M_x$. It is independent of the choice of the basis.

\section{Matrix Traces and the Radical}\label{mat}
Our main construction is based on the following results describing
the elements of the radical of an associative algebra $A$ using
traces of elements:

\begin{thm}[Dickson \cite{dickson} pp.106-107] \label{dic}
An element $x$ of an associative algebra $A$ over a field $F$ of
characteristic $0$ is properly nilpotent if and only if $Tr(xy)=0$
for every $y\in A$.
\end{thm}

\begin{cor}[Friedl and R\'onyai \cite{ronyai85} p.156] \label{ron}
Let $F$ be a field of characteristic $0$ and $A$ a matrix algebra
over $F$. Let $B=[b_1,\ldots,b_n]$ be a linear basis of $A$ over
the field $F$. Then $x\in Rad(A)$ if and only if $Tr(xb_i)=0,\quad
i=1,\ldots,n.$
\end{cor}

We apply the above results to the special case of commutative
algebras which are quotients of polynomial rings. Consider the
system of polynomial equations
\[
\bold{f}(\bold{x})=0
\]
where $\bold{f}=\{f_1,\ldots,f_l\}$ and each $f_i$ is a polynomial
in the variables $\bold{x}=[x_1,\ldots,x_m]$. Assume that the
polynomials $f_1,\ldots,f_l$ have finitely many roots in
$\mathbb{C}^m$, which implies that the algebra
$A=\mathbb{C}[\bold{x}]/I$ is finite dimensional, where $I$ is the
ideal generated by the polynomials in $\bold{f}$. Denote the
dimension of $A$ over $\mathbb C$ by $n$ and let
$B=[b_1,\ldots,b_n]$ be a basis of $A$. By slight abuse of
notation we denote the elements of the basis $B$ which are in $A$
and some fixed preimages of them in $\mathbb{C}[x_1, \ldots, x_m]$
both by $b_1, \ldots, b_n$. Let
$\{\bold{z}_1,\ldots,\bold{z}_n\}\subset \mathbb{C}^m$ be the set
of common roots (not necessarily all distinct) of the polynomials
in $\bold{f}$. Using the multiplication matrices $M_f$ associated
to the elements $f\in A$ and the fact that
\[
Rad(A)=\sqrt{I}/I\subseteq\mathbb{C}[\bold{x}]/I=A,
\]
we can reword Corollary \ref{ron} in the following way:
\begin{cor}\label{rad}
Let $p\in \mathbb{C}[\bold{x}]$ and $\bar p$ be the image of $p$ in
$A$. Using the above notation, the following statements are
equivalent:

\begin{enumerate}
\item[(\it{i})] $p\in \sqrt{I}$
\item[(\it{ii})]$\bar p\in Rad(A)$
\item[(\it{iii})] $Tr(M_{\bar{p}b_j})=0$ for all $ j=1,\ldots,n.$
\end{enumerate}
\end{cor}

We can now use the previous corollary to characterize the radical
of $A$ as the nullspace of a matrix defined as follows:

\begin{defn}\label{tracemat}
The \emph{matrix of traces} is the $n\times n$ symmetric matrix:
\[
R=\left[Tr(M_{b_ib_j})\right]_{i,j=1}^n
\]
where $M_{b_ib_j}$ is the multiplication matrix of $b_ib_j$ as an
element in $A$ in terms of the basis $B=[b_1, \ldots, b_n]$ and
$Tr$ indicates the trace of a matrix.
\end{defn}

\begin{cor}\label{sol}
An element
\[
r=\sum_{k=1}^{n} c_k b_k
\]
of the quotient ring $A$ with basis $B=[b_1,\ldots,b_n]$ is in the
radical of $A$ if and only if $[c_1,\ldots,c_n]$ is in the
nullspace of the matrix of traces $R$.
\end{cor}

\begin{proof}
Corollary \ref{rad} states that an element $r=\sum_{k=1}^{n} c_k
b_k\in A$ belongs to $Rad(A)$ if and only if $Tr(M_{rb_j})=0$, for
all $j=1,\ldots,n.$ From the linearity of both the multiplication
map (see Proposition (4.2) in Chapter 2 of \cite{Cox98}) and the
traces of matrices we have that
\[
\begin{aligned}
Tr(M_{rb_j})&=\sum_{k=1}^{n} c_k Tr(M_{b_kb_j})\notag\\
            &=[c_1,\ldots,c_n]R[j]\notag\\
\end{aligned}
\]
where $R[j]$ is the $j^{th}$ column of the matrix of traces $R$.
Therefore, $Tr(M_{rb_j})=0$ for all $j=1,\ldots,n$ is equivalent to
$[c_1,\ldots,c_n]R=0$.

\end{proof}

\begin{rem}\label{rmmulti}\emph{
Methods in the literature for computing the matrix of traces $R$
are mentioned in the Introduction. One way to compute it is from
the multiplication matrices $M_{b_ib_j}$. Note that in order to
compute the matrices $M_{b_ib_j},\; i,j=1,\ldots,n$, it is
sufficient to have $M_{x_k},\; k=1,\ldots,m$, since if $h\in
\mathbb{C}[x_1, \ldots, x_m]$  is a preimage of $b_ib_j\in A$,
then we have
\[
\begin{aligned}
 M_{b_ib_j}&=M_{h(x_1,\ldots,x_m)}\notag\\
           &=h(M_{x_1},\ldots,M_{x_m}).\notag\\
\end{aligned}
\]
This is because the regular representation is a homomorphism of
$\mathbb C$-algebras, see also Corollary (4.3) in Chapter 2 of
\cite{Cox98}.}
\end{rem}

\begin{ex}\label{ex1}
We consider the polynomial system $f_1=f_2=f_3=0$, with
\begin{scriptsize}
\[
\begin{aligned}
f_1=&x_1^2+4x_1x_2-6x_1+6x_2^2-18x_2+13\notag\\
f_2=&x_1^3+16x_1^2x_2-7x_1^2+118x_1x_2^2-286x_1x_2\notag\\&+147x_1-x_2^3+6x_2^2+x_2+5\notag\\
f_3=&x_1^3+10x_1^2x_2-5x_1^2+72x_1x_2^2-176x_1x_2\notag\\&+91x_1-x_2^3+4x_2^2+x_2+3\notag
\end{aligned}
\]
\end{scriptsize}
These polynomials have two common roots: $[1,1]$ of multiplicity 3
and $[-1,2]$ of multiplicity 2.

We compute the multiplication matrices $M_{x_1}$ and $M_{x_2}$
with respect to the basis $B= [1,x_1,x_2,x_1x_2,x_1^{2}]$, which
are respectively
\begin{tiny}
\begin{eqnarray}\label{multi1}
\left[ \begin {array}{rrrrr} 0&1&0&0&0\\\noalign{\medskip}
0&0&0&0&1\\\noalign{\medskip} 0&0&0&1&0\\\noalign{\medskip}
\frac{5}{3}&-2&-1&\frac{2}{3}&\frac{5}{3}\\\noalign{\medskip}
-{\frac {17}{3}}&1&4&\frac{4}{3}&\frac{1}{3}\end {array} \right], \,
 \left[
\begin {array}{rrrrr}
0&0&1&0&0\\\noalign{\medskip}0&0&0&1&0\\\noalign{\medskip}-{\frac
{13}{6}}&1&3&-{\frac{2}{3}}&-{\frac{1}{6}}\\\noalign{\medskip}-{\frac{1}{6}}&-1&0&\frac{7}{3}&-{\frac{1}{6}}\\\noalign{\medskip}\frac{5}{3}&-2&-1&\frac{2}{3}&\frac{5}{3}\end
{array} \right].
\end{eqnarray}
\end{tiny}

Here we used Chardin's subresultant construction to compute the
multiplication matrices. (See \cite{Cha} and \cite{Sza2}.)

We now compute the matrix $R$ using Definition \ref{tracemat} and
Remark \ref{rmmulti}:
\begin{scriptsize}
\begin{eqnarray}\label{rmat}
R= \, \left[ \begin {array}{rrrrr}
5&1&7&-1&5\\\noalign{\medskip}1&5&-1&7&1\\\noalign{\medskip}7&-1&11&-5&7\\\noalign{\medskip}-1&7&-5&11&-1\\\noalign{\medskip}5&1&7&-1&5\end
{array} \right].
\end{eqnarray}
\end{scriptsize}

The nullspace of R is generated by the vectors
\begin{scriptsize}$$[1, -3, 0, 2, 0], [0, -4, 1, 3, 0], [0, -3, 0, 2,
1].$$\end{scriptsize} By Corollary \ref{sol} we have that the
radical  of $I=\left<f_1,f_2,f_3\right>$ modulo $I$ is
\begin{scriptsize}
\[
\sqrt{I}/I=\left<1-3x_1+2x_1x_2, -4x_1+x_2+3x_1x_2,
-3x_1+2x_1x_2+x_1^2\right>.
\]
\end{scriptsize}
Note that the polynomials on the right hand side are in
$\sqrt{I}$.
\end{ex}

Assume that ${\rm rank} \;R=k$. Once we know the $n-k$ generators
$\{r_{k+1},\ldots,r_n\}$ of the radical, we can obtain the
multiplication matrices of the elements of
$A/Rad(A)=\mathbb{C}[\bold{x}]/\sqrt{I}$ by performing a change of
basis on the multiplication matrices $M_{x_1},\ldots,M_{x_m}$ to
the basis $\{r_1,\ldots,r_k,r_{k+1},\ldots,r_n\}$ of $A$, where
$r_{1},\ldots,r_k$ can be chosen arbitrarily as long as
$\{r_1,\ldots,r_k,r_{k+1},\ldots,r_n\}$ is linearly independent.
Let $M_{x_s}$ be the multiplication matrix of the coordinate $x_s$
in the basis $[r_1,\ldots,r_n]$. Then the $k\times k$ principal
submatrix
\[
M'_{x_s}:=\left[M_{x_s}(i,j)\right]_{i,j=1}^k
\]
is the multiplication matrix of $x_s$ in
$A/Rad(A)=\mathbb{C}[\bold{x}]/\sqrt{I}$ with respect to the basis
$[r_{1},\ldots,r_k]$.

\begin{ex}\label{ex2}
Continuing Example \ref{ex1}, we have that the generators of the
radical $Rad(A)$ have coordinates
\begin{scriptsize}
\[
r_{3}=[1, -3, 0, 2, 0],\; r_{4}=[0, -4, 1, 3, 0],\; r_5=[0, -3, 0,
2, 1]
\]
\end{scriptsize}
in the basis $B=[1,x_1,x_2,x_1x_2,x_1^{2}]$.

We set \begin{scriptsize}
\[
r_1=[1, 0, 0, 0, 0],\; r_2=[0, 1, 0, 0, 0].
\]
\end{scriptsize}

We perform the change of basis to the two multiplication matrices
$M_{x_1}$ and $M_{x_2}$ and obtain:
\begin{footnotesize}
\[
\left[\begin{array}{rrrrr}
 0& 1&0& 0& 0\\
1& 0&-1& 0& 1\\
0&0& 10/3& -2& 1/3\\
 0&0&5& -3&1\\
 0& 0&-7/3& 2& 2/3
\end{array} \right]
\; \text{ and }\; \left[ \begin {array}{rrrrr}
 3/2& -1/2&-3/2& 1& 0\\
 -1/2& 3/2&1/2& 0& 0\\
  0& 0&-1/3& 1& -1/3\\
 0& 0&-8/3& 3& -2/3\\
 0& 0&4/3& -1& 4/3
\end {array} \right]
\]
\end{footnotesize}
respectively.

We then have that the multiplication matrices for $x_1$ and $x_2$ in
$A/Rad(A)$ in the basis $[1,x_1]$ are
\begin{scriptsize}
\[
\mathcal{M}_{x_1}=\left[\begin{array}{rr}0&1\\1&0\end{array}\right]
\; \text{ and }\;
\mathcal{M}_{x_2}=\left[\begin{array}{rr}3/2&-1/2\\-1/2&3/2\end{array}\right].
\]
\end{scriptsize}
The eigenvalues of these matrices give the solutions to the system.
\end{ex}

\section{Clustered roots}

In this section we  consider systems with clustered roots instead of
systems with root multiplicities. We can think of these systems with
clustered roots as being obtained from systems with multiplicities
via one of the following two ways:
\begin{enumerate}
\item by perturbing the coefficients of the system with multiple
roots, \item by perturbing the multiple roots to obtain clusters.
\end{enumerate}

Let $\bold{f}$ be the system with multiple roots and
$\tilde{\bold{f}}$ be the system with clustered roots obtained from
$\bold{f}$ by any of the above methods. Denote by
$\tilde{A}=\mathbb{C}[\bold{x}]/\tilde{I}$ the algebra corresponding
to the ideal $\tilde I$ generated by the polynomials in $\tilde{\bold{f}}$.\\

\begin{assu} Throughout this paper we make the
assumption that the basis $B$ for $A$ also forms a basis for
$\tilde{A}$. Note that if $\bold{f}$ is a well constrained system
then for ``almost all" perturbations $\tilde{\bold{f}}$ our
assumption is satisfied, i.e. the set of perturbed systems for
which it doesn't hold has measure zero in the space of all systems
of given degrees.\end{assu}

If we assume that the basis $B$ for $A$ also forms a basis for
$\tilde{A}$ then both the multiplication matrices and the matrix of
traces are continuous functions of the coefficients of the
polynomials. Therefore, small perturbations in the coefficients of
$\bold{f}$ will result in small changes in the entries of the
multiplication matrices and the matrix of traces.

However, in case 2, when the roots are perturbed,  the polynomials
corresponding to the clustered system might end up having
coefficients very different to those of the original system, even
if the radii of the clusters were small. In this case, if we
compute the multiplication matrices for the clustered system, the
entries might not be continuous functions of the perturbation of
the roots. They not only depend on the magnitude of the
perturbation of the roots but also on the \textit{direction} of
the perturbation. However, as we shall show in  Proposition
\ref{traceroots}, the matrix of traces is always continuous in the
roots. The following examples illustrates this phenomenon.

\begin{ex}\label{clusterex}
We consider three examples of a single cluster of size
proportional to $\varepsilon$ around the origin $(0,0)$ in
$\mathbb{C}^2$ consisting of three roots. The first two examples
demonstrate that the defining equations and the multiplication
matrices can have different accumulation points as $\varepsilon$
approaches $0$, depending on the direction. The third example
demonstrate that generally the defining equations and the
multiplication matrices are not continuous at $\varepsilon=0$.
 \begin{itemize}

 \item First, the roots of the cluster are
 $
 (0,0), \;(\varepsilon, \varepsilon),
 \;(2\varepsilon,2\varepsilon).
 $ The defining equations of these points in $\mathbb{C}[x, y]$ are given
 by $x^3-3\varepsilon x^2+2\varepsilon^2x=0$ and $y=x$, and the
 multiplication matrices in the basis $B=\{1,x,x^2\}$ are given by
 $$
 M_x=M_y=\left[\begin{array}{rrr}0&1&0\\0&0&1\\0 &-2\varepsilon^2
 &3\varepsilon\end{array}\right] \quad \text{ and } \quad
 \lim_{\varepsilon\rightarrow 0} M_y=\left[\begin{array}{rrr}0&1&0\\0&0&1\\0
 &0
 &0\end{array}\right],
 $$
 and the primary ideal defining the multiple root is $\langle x^3,
 x-y\rangle$.

 \item The next example has cluster $
 (0,0), \;(\varepsilon, 2\varepsilon),
 \;(2\varepsilon,4\varepsilon).
 $ The defining equations are $x^3-3\varepsilon x^2+2\varepsilon^2x=0$ and $y=2x$. Then $M_x$ is the same as above, but
 $$
 M_y=\left[\begin{array}{rrr}0&2&0\\0&0&2\\0 &-4\varepsilon^2
 &6\varepsilon\end{array}\right] \quad \text{ and } \quad
 \lim_{\varepsilon\rightarrow 0} M_y=\left[\begin{array}{rrr}0&2&0\\0&0&2\\0
 &0
 &0\end{array}\right],
 $$
 and the primary ideal defining the multiple root is $\langle x^3,
 x-2y\rangle$.

\item More generally, the third example has cluster $
 (0,0), \;(\varepsilon, c\varepsilon),
 \;(2\varepsilon,d\varepsilon)$ for some $c,d\in \mathbb{R}$.
 Then the first defining equation is the same as above, and the second
 equation is $y= -\frac{2c-d}{2\varepsilon}x^2 + (2c-\frac{d}{2})x$ which is not
 continuous in $\varepsilon=0$, unless $d=2c$. Similarly for the multiplication matrix $M_y$.  However, the matrix of traces
 $$
 R= \left[\begin{array}{rrr}3&c\varepsilon+d\varepsilon&c^2\varepsilon^2+d^2\varepsilon^2\\c\varepsilon+d\varepsilon&
 c^2\varepsilon^2+d^2\varepsilon^2&c^3\varepsilon^3+d^3\varepsilon^3\\c^2\varepsilon^2+d^2\varepsilon^2 &c^3\varepsilon^3+d^3\varepsilon^3
  &c^4\varepsilon^4+d^4\varepsilon^4\end{array}\right]
  $$ (with respect to the basis $\{1,y,y^2\}$) is continuous in $\varepsilon=0$ and has the same limit
   for every choices of $c$ and
  $d$.
 \end{itemize}

\end{ex}
\begin{ex}\label{ex3} Continuing with Example \ref{ex2},
suppose now that instead of having a system with common roots
$[1,1]$ of multiplicity 3 and $[-1,2]$ of multiplicity 2 we have a
polynomial system with a cluster of three common roots:
\[
[[1, 1], [0.9924, 1.0027], [1.0076, 0.9973]]
\]
around $[1,1]$ and a cluster of two common roots:
\[
[[-1, 2], [-1.0076, 2.0027]]
\]
around $[-1,2]$.

Using the multivariate Vandermonde construction (see for example
\cite{moste95}), we obtained the following multiplication matrices
for this system, with respect to the same basis as for the system
with multiple roots: $B= [1,x_1,x_2,x_1x_2,x_1^{2}]$.
\begin{scriptsize}
\[
\begin{aligned}
\tilde{M}_{x_1}&= \, \left[
\begin {array}{rrrrr}
3.8328\times 10^{-6}& 9.9997\times 10^{-1}&
3.0830\times 10^{-8} &4.1421\times 10^{-7}& 3.9951\times 10^{-5}\\
3.7919\times 10^{-6}&-2.7338\times 10^{-5}&
1.2303\times 10^{-7} &8.2891\times 10^{-7}& 1.00004\\
3.8527\times 10^{-6}&-
2.7463\times 10^{-5}&- 1.5183\times 10^{-8} &1.00000& 3.9969\times 10^{-5}\\
7.73947& 21.69983&- 5.97279 &- 16.79084&- 5.67565\\
- 17.94136&- 54.43008& 13.97610 &41.61207& 17.78328
\end {array}
\right]\notag\\
\end{aligned}
\]
\[
\begin{aligned}
\tilde{M}_{x_2}&= \, \left[
\begin {array}{rrrrr}
3.7831\times 10^{-6}&- 2.7103\times 10^{-5}&
1.00000 &- 1.4715\times 10^{-7}& 4.0017\times 10^{-5} \\
3.8527\times 10^{-6}&- 2.7464\times 10^{-5}&-
1.5183\times 10^{-8} &1.00000& 3.9969\times 10^{-5}\\
- 2.22905&
1.06576& 3.00000 &- 7.1053\times 10^{-1}&- 1.2617\times 10^{-1}\\
- 3.23468&- 10.77768&
2.47988 &9.67839& 2.85410\\
7.73947& 21.69983&- 5.97279 &- 16.79084&- 5.67565
\end {array}
\right]\notag\\
\end{aligned}
\]
\end{scriptsize}

The norm of the difference between these matrices and the
multiplication matrices (\ref{multi1}) for the system with multiple
roots are very large: $135.41$ for the matrices of $x_1$ and $59.54$
for the matrices of $x_2$. Entrywise, the largest absolute value of
the difference of the entries of the matrices is $55.40$ for $x_1$
and $21.70$ for $x_2$.

However, the matrix of traces associated to the system with clusters
is
\begin{scriptsize}
\begin{equation}\label{rmatc}
\left[ \begin {array}{rrrrr}
4.99999& 0.99240&7.00269&- 1.01796& 5.01538\\
0.99259& 5.01557&- 1.01777& 7.03349& 0.97757\\
7.00131&- 1.01934& 11.00943&- 5.04274& 7.03192\\
-1.01900& 7.03226&- 5.04240& 11.07093&- 1.04951\\
5.01548& 0.97748& 7.03339&- 1.04838& 5.03155
\end{array} \right]
\end{equation}
\end{scriptsize}
and the 2-norm of the difference between this matrix and the
multiplication matrix $R$ in (\ref{rmat}) for the system with
multiple roots is $0.147$.
\end{ex}

We have the following result for the entries of the matrix of traces
$R$ expressed in terms of the roots of the polynomial system.

\begin{prop}\label{traceroots}
The matrix of traces $R$ of the system $\bold{f}(\bold{x})=0$ with
respect to $B=[b_1, \ldots, b_n]$ can be expressed in terms of the
common roots $\{\bold{z}_1,\ldots,\bold{z}_n\}$   as
\[
R=\left[\sum_{k=1}^n b_ib_j(\bold{z}_k)\right]_{i,j=1}^n
\]
where $b_ib_j(\bold{z}_k)$ indicates the evaluation of the
polynomial $b_ib_j$ at the point $\bold{z}_k$.
\end{prop}

\begin{proof}
Assume that $\{\bold{z}_1,\ldots,\bold{z}_k\}$ are the distinct
elements among $\{\bold{z}_1,\ldots,\bold{z}_n\}$ in $V(I)$ and let
$n_i$ be the multiplicity of $\bold{z}_i$. Let $Q_i$ be the (unique)
primary component of $I$ in ${\mathbb C}[\bold{x}]$ whose radical
$P_i$ is the ideal of all polynomials vanishing at $\bold{z}_i$,
$i=1,2,\ldots ,k$. Set $A_i={\mathbb C}[\bold{x}]/Q_i$. We have then
$n_i=\dim _{\mathbb C} A_i$ and $I=Q_1\cap Q_2\cap\cdots Q_k$. Also
the ideals $Q_i$ are pairwise relatively prime, hence by the Chinese
Remainder Theorem we have
\[
A\cong A_1\oplus A_2\oplus\cdots \oplus A_k.
\]
We denote also by $A_i$ the image of $A_i$ in $A$ at this
isomorphism. Given any polynomial $g$, it is immediate that $A_i$
is an invariant subspace of the multiplication map $M_g$ and that
the characteristic polynomial of $M_g$ on $A_i$ is
$\left(t-g(\bold{z}_i)\right)^{n_i}$. This implies that the
characteristic polynomial of $M_g$ is $\prod_{i=1}^k
\left(t-g(\bold{z}_i)\right)^{n_i}=$ $\prod_{i=1}^n
\left(t-g(\bold{z}_i)\right)$. So the trace of $M_g$ is
$\sum_{i=1}^ng(\bold{z}_i)$. Therefore

\[
Tr(M_{b_ib_j})=\sum_{k=1}^n b_ib_j(\bold{z}_k)
\]
which proves the lemma.
\end{proof}

Note: An alternative proof can be given for Proposition
\ref{traceroots} using the fact that the multiplication matrix $M_g$
is similar to a block diagonal matrix where the $i$-th diagonal
block is an $n_i\times n_i$ upper triangular matrix, with diagonal
entries $g(\bold{z}_i),\;i=1,\ldots,k$ (cf. \cite[Theorem
2]{moste95}).

The previous result shows that the entries of the matrix of traces
are continuous functions of the roots, even when the roots
coincide. In particular, a system with multiple roots and a system
with clusters obtained by perturbing the roots of  a system with
multiplicities will have comparable matrices of traces.

\section{Univariate Case}
Before we give our method in full generality we would like to
describe our algorithm in the univariate case. The purpose of this
section is to demonstrate the simplicity and the accuracy of our
technique to compute the approximate square-free factorization of a
univariate polynomial. As we mentioned in the Introduction, our
method offers a new alternative to other approximate square-free
factorization algorithms, such as the one in \cite{KalMay2003}.

The following is a description of the steps of our algorithm. Let
\[
f(x)=x^d+a_{1}x^{d-1}+\cdots+a_{d-1}x+a_d\in \mathbb{C}[x]
\] be a given polynomial of degree $d$ with clusters of roots of size at
most $\varepsilon$. The output of our algorithm is a polynomial
$g(x)\in \mathbb{C}[x]$ such that its roots are the arithmetic
means of the roots in each cluster, with a precision of order of
magnitude $\varepsilon^2$.
\begin{enumerate}
\item  Compute the matrix of traces $R$ w.r.t. the basis
$B=[1,x,x^2,\ldots,x^{d-1}]$ using the Newton-Girard formulas. In
this case we have $R=\left[s_{i+j}\right]_{i,j=0}^{d-1}$ where
$s_{t}$ is the sum of the $t$-th power of the roots of $f$. We set
$s_0=d$ and we find $s_1, \ldots, s_{2d-2}$ from the coefficients
of $f$ using the Newton-Girard formulas as follows:
\begin{footnotesize}
\[
\begin{aligned}
s_1+a_1&=0\notag\\
s_2+a_{1}s_1+2a_{2}&=0\notag\\
\vdots &\notag\\
s_d+a_{1}s_{d-1}+\cdots+a_{d-1}s_1+da_d&=0\notag\\
s_{d+1}+ a_1s_{d}+\cdots+a_{d}s_1&=0\notag\\
\vdots &\notag\\
s_{2d-2}+a_1s_{2d-3}+\cdots+a_{d}s_{d-3}&=0.\notag
\end{aligned}
\]
\end{footnotesize}

\item  Gaussian elimination with complete pivoting (GECP) is used
on the matrix $R$ until the remaining entries in the partially row
reduced matrix $U_k$ are smaller than a preset threshold (see
Propositions \ref{elemrow} and \ref{bound}). The number of
iterations performed, $k$, is the numerical rank of the matrix $R$.

\item Compute a basis  of the nullspace  $N$ of the first $k$ rows
of the matrix $U_k$ obtained after $k$ steps of the GECP. We
identify the vectors in $N$ by polynomials,  by combining their
coordinates with the corresponding basis elements of $B$.

 \item The smallest degree polynomial
in $N$ is the approximate square-free factor $g(x)$ of $f(x)$. Its
roots are the arithmetic means of the roots in each cluster modulo
$\varepsilon^2$ (see Proposition \ref{meangen}). In the case when
the matrix $R$ has numerical rank $d$ then we take $g(x)=f(x)$
 as the square-free factor.
\end{enumerate}

\begin{ex}
(1) Consider the approximate polynomial
\begin{footnotesize}
\[
f=(x-(z+\delta_1\varepsilon))(x-(z+\delta_2\varepsilon))(x-(z+\delta_3\varepsilon))
\]
\end{footnotesize}
obtained by perturbing the roots of the polynomial
\[
(x-z)^3=x^3-3x^2z+3xz^2-z^3.
\]

Using the basis $B=[1,x,x^2]$ we obtained the matrix of traces
$R$, for which the $U$ matrix in the $LU$ factorization obtained
by GECP is
\begin{tiny}
\[
\left(
\begin{array}{ccc}
3& 3z+\varepsilon\left(\delta_3+\delta_2 +\delta_1\right)
&3z^2+\varepsilon\left(2z\delta_3+2\delta_2 z+2\delta_1
z+\right)+\varepsilon^2\left(\delta_3^2+\delta_2^2+\delta_1^2\right)\\
0&
\frac{\varepsilon^2\left(-2\delta_1\delta_3-2\delta_2\delta_3-2\delta_1\delta_2+2\delta_1^2+2\delta_2^2+2\delta_3^2\right)}{3}
&\frac{\varepsilon^2\Phi_{2,2}+\varepsilon^3\Phi_{2,3}}{3}\\
0& 0&\frac{\varepsilon
^4\Phi_{3,3}}{2\left(-\delta_1\delta_3-\delta_2\delta_3-\delta_1\delta_2+\delta_1^2+\delta_2^2+\delta_3^2\right)}
\end{array}
 \right)
\]
\end{tiny}
where $\Phi_{i,j}$ are polynomials in the $\delta$'s and $z$'s.

Using the bound from Proposition \ref{bound} for the numerical rank,
we have that the approximate radical will be defined using the
nullspace of the first row of $R$.

We obtain the following basis of the approximate radical,
\begin{small}
\[
\{x^2-z^2-\frac{2z\varepsilon\left(\delta_3
+\delta_2+\delta_1\right)+\varepsilon^2\left(\delta_3^2+\delta_2^2+\delta_1^2\right)}{3},
x-z-\frac{\varepsilon\left(\delta_3+\delta_2 +\delta_1\right)}{3}\}
\]
\end{small}

We choose the element of smallest degree to be the approximate
square-free factor of $f$, which is here
\begin{small}
\[
x-z-\frac{\varepsilon\left(\delta_3+\delta_2 +\delta_1\right)}{3}.
\]
\end{small}

We can see that in this case the roots of this polynomial correspond
precisely to the arithmetic mean of the three clustered roots.
\\

 (2) Consider the approximate polynomial
\begin{scriptsize}
\[
\begin{aligned}
f(x)&=
(x+(-0.98816+0.01847I))(x+(-0.98816-0.01847I))\notag\\
&\qquad(x-1.02390)(x-1.98603)(x-2.01375)\notag
\end{aligned}
\]
\end{scriptsize}
which is a perturbation of the polynomial
\begin{scriptsize}
\[
x^5-7x^4+19x^3-25x^2+16x-4=(x-1)^3(x-2)^2.
\]
\end{scriptsize}

The matrix of traces corresponding to $f$ is
\begin{tiny}
\[
R=\left[\begin{array}{rrrrrr}5& 7.00001& 11.00013&
19.00089& 35.00425\\
7.00001& 11.00013& 19.00089& 35.00425&
67.01631\\
11.00013& 19.00089& 35.00425& 67.01631&
131.05456\\
19.00089& 35.00425& 67.01631& 131.05456&
259.16598\\
35.00425& 67.01631& 131.05456& 259.16598& 515.47172
\end{array}\right].
\]
\end{tiny}

The $U_k$ matrix obtained after 2 steps of GECP on $R$ is
\begin{tiny}
\[
U_2=\left[\begin{array}{rrrrrr} 515.47172& 35.00425& 131.05456&
259.16598& 67.01631\\
0& 2.62296& 2.10058& 1.40165& 2.44912\\
0& 0& 0.0024342& 0.0029279& 0.0011698\\
0& 0& 0.0029279& 0.0035326& 0.0014044\\
0& 0& 0.0011698& 0.0014044& 0.00056307
\end{array}\right].
\]
\end{tiny}

By taking the nullspace of the first two rows of the matrix $U_2$,
we obtain the following basis of the approximate radical,
\begin{scriptsize}
\[\begin{aligned}
&\{x^4-15.01431x+14.01921, x^3-7.00397x+6.00539,\\
&\quad x^2-3.00074x+2.00102\}.\\
\end{aligned}\]
\end{scriptsize}

The approximate square-free factor of $f$ is then
\begin{scriptsize}
\[
x^2-3.00074x+2.00102=(x-1.00028)(x-2.00047).
\]
\end{scriptsize}
We can see that the roots of the output are close to the means of
the clusters, and the differences are $0.00058$ and $0.000200$
respectively, which are  of the order of the square of the cluster
size (bounded here by $0.03$).
\end{ex}

We refer to the papers of
\cite{SasNod89,hriste97,KalMay2003,Zeng2003}  for other methods
that study approximate square-free factorization using approximate
gcd computation.

\section{LU decomposition of the matrix of traces}

Since the polynomial system with clusters, obtained by perturbing
the system with multiplicities, has only simple roots, the matrix
of traces has full rank. However, we can try to find its numerical
rank. We will argue below that we can define the numerical rank in
such a way that it will be equal to the rank of the matrix of
traces of the corresponding system with multiplicities.

In this paper we primarily study the Gaussian elimination with
complete pivoting (GECP) \cite{GoVa96} in order to estimate the
numerical rank and find the numerical nullspace of the matrix of
traces. However we we will also infer that the singular value
decomposition (SVD) in our case works similarly to the GECP.

We would like to note that rounding errors can sometimes result in
a matrix which is close to a singular one, but where all the
pivots are large (see Kahan's Example 5 in \cite{Kahan66}). This
example shows that GECP can be a poor choice for numerical rank
computations in the presence of rounding errors. On the other
hand, algorithms for the accurate computations of the SVD of
certain structured matrices, including Vandermonde matrices, use
improved versions of GECP as subroutines \cite{Demmel99,DeKoe05}.
In our case we prove that the structure of the matrix of traces
guarantees that we will obtain small pivots which are proportional
to the square of the size of the clusters and can therefore use
GECP for rank determination.

We will also show how our results for the GECP of the matrix of
traces $R$ relate to the singular values of $R$. In particular we
will obtain asymptotic bounds for  the singular values of the
matrix $R$. Such bounds are similar to the ones for the entries of
the $U_k$ matrix obtained after $k$ steps of GECP on $R$, more
precisely, we also obtain in this case that the ``almost zero"
singular values are proportional to the square of the size of the
clusters.

First we study the properties of the Gaussian elimination in the
approximate setting. We use the following notation for different
versions of the Gaussian elimination algorithm:

\begin{defn}\label{Gauss} The version of Gaussian elimination in which at the
$i$-th step we always select the entry at position $(i,i)$ for
pivoting will be referred to as {\em regular}. We call an $m\times
n$ matrix $M$ {\em regular} if for $k := {\rm rank}(M)$ the first
$k$ steps of the regular Gaussian elimination on $M$ do not
 encounter zero pivots.

Note that GECP on the matrix $M$ computes two permutation matrices
$P$ and $Q$ of sizes $m\times m$ and $n\times n$, respectively,
such that for the matrix $P\,M\,Q$ the regular Gaussian
elimination works as GECP.
 \end{defn}

In the rest of this section we give results which compare the GECP
applied to  the matrices of traces of the perturbed system and to
 the system with multiple roots. Let $R_0$ be the matrix of
traces of the system with multiple roots and let  ${R}$ denote the
matrix of traces of some perturbation of it. Assume that ${\rm
rank}(R_0)=k$. Our next result guarantees that for sufficiently
small clusters, the first $k$ steps of the  GECP applied to ${R}$
computes permutation matrices $P$ and $Q$ which make the matrix
$P\,R_0\,Q$ regular.

\begin{prop}\label{regular}
Let $M$ be an $n\times n$ matrix with entries polynomials in ${\bf
x}=[x_1, \ldots, x_N]$ over $\mathbb{C}$. Fix ${\bf z}=[z_1,
\ldots, z_N]\in \mathbb{C}^N$, denote $M_0:= M|_{{\bf x}={\bf
z}}$, and assume that ${\rm rank}(M_0)=k$. Then there exists an
open neighborhood ${\mathcal V}$ of ${\mathbf z}$ in
$\mathbb{C}^N$ such that for all points $\tilde{\bf
z}=[\tilde{z}_1, \ldots, \tilde{z}_N]\in {\mathcal V}$ if $P$ and
$Q$ are the permutation matrices corresponding to the first $k$
steps of the GECP on the matrix $\tilde{M}:=M|_{{\bf x}=\tilde{\bf
z}}$, then
   the matrix $P\;M_{0}\;Q$ is
regular.
\end{prop}

\begin{proof}
We call a pair $(P,Q)$ of $n$ by $n$ permutation matrices {\em
good} if $P\,M_{0}\,Q $ is regular, otherwise the pair is called
{\em bad}. For each bad pair we define an open neighborhood
${\mathcal V}_{P,Q}$ of ${\mathbf z}\in \mathbb{C}^N$ as follows:
  For some $i\leq k$
assume that the regular Gaussian  elimination on  $P\, M_{0}\, Q $
encounters a  zero pivot for the first time in the $i$-th step,
causing $(P,Q)$ to be a bad pair. Denote by $U_0$ the partially
reduced form of $P\,M_{0}\,Q $ after the $i-1$-th step of the
regular Gaussian elimination. Denote by $S$ the set of indices
$(s,t)$ such that $s,t\ge i$ and the $(s, t)$ entry of $U_0$ is
non-zero, and by $T$ the set of indices $(s,t)$ such that $s,t\ge i$
and the $(s, t)$ entry of $U_0$ is zero. Since the rank of
$P\,M_{0}\,Q $ is $k$, $S$ is non empty.

Let $U$ be the partially reduced matrix obtained from $P\,M\,Q $ via
the first $i-1$ steps of regular Gaussian  elimination. Note that
the entries of $U$ are rational functions of the entries of $M$  and
the denominators of these are non zero at ${\bf z}$, hence  are
continuous functions of the  points $[\tilde{z}_1, \ldots,
\tilde{z}_N]$ in a sufficiently small neighborhood of ${\mathbf z}$.
In particular, in an open neighborhood ${\mathcal U}$ of ${\mathbf
z}$ the first $i-1$ steps of regular elimination can be carried out.

Let the open neighborhood
  ${\mathcal V}_{P,Q}\subset {\mathcal U}\subset \mathbb{C}^N$ of
${\mathbf z}$ be selected such that for all $[\tilde{z}_1, \ldots,
\tilde{z}_N] \in {\mathcal V}_{P,Q}$  the entries in $T$ of
$\tilde{U}:=U|_{{\bf x}=\tilde{\bf z}}$ are all strictly smaller in
absolute value than any of the entries in $S$ of $\tilde{U}$.
 By continuity, such open neighborhood of
${\mathbf z}$ exists, since the required inequalities hold for
$U_0$.

Finally define $\displaystyle{{\mathcal V}:=\bigcap_{(P,Q)\;
\text{ is bad}} {\mathcal V}_{P,Q}}$. This is also an open
neighborhood of ${\mathbf z}$ since the set of permutations is
finite. We claim that for any fixed $[\tilde{z}_1, \ldots,
\tilde{z}_N] \in {\mathcal V}$, if $(P, Q)$ is the pair of
permutation matrices corresponding to the first $k$ steps of the
Gaussian elimination with complete pivoting on the matrix
$\tilde{M}$ then $(P,Q)$ is a good pair. This is true since
$P\tilde{M}Q$ has the property that for each $i\leq k$ after $i-1$
steps of the Gauss elimination the $(i,i)$-th  entry of the
corresponding matrix is maximal in absolute value among the
entries indexed by $(s,t)\neq (i,i)$ such that $s,t\geq i$. But
then the $(i,i)$-th entry in the matrix $U_0$ defined above cannot
be $0$ because of the definition of ${\mathcal V}$. This proves
the claim.
\end{proof}

In the rest of the paper we will assume that the size of the
clusters is a parameter   $\varepsilon$. More precisely, in the
following definition we formally explain the mathematical setting
where our results will hold:

\begin{defn} \label{setting}  Let
$\bold{z}_i=[z_{i,1},\ldots,z_{i,m}]\in \mathbb{C}^m$ for
$i=1,\ldots,k$, and consider $k$ clusters  $C_1, \ldots, C_k$  of
size $|C_i|=n_i$ such that $\sum_{i=1}^kn_i=n$, each of radius
proportional to the parameter $\varepsilon$ in the $\infty$-norm
around $\bold{z}_1, \ldots, \bold{z}_k$:
\begin{equation}\label{cluster}
\begin{scriptsize}
\begin{aligned}
C_i=&\{[z_{i,1}+\delta_{i,1,1}\varepsilon,\ldots,z_{i,m}+\delta_{i,
1,m}\varepsilon],
\ldots,\\
&\ldots[z_{i,1}+\delta_{i, n_i,1}\varepsilon,\ldots,z_{i,m}+\delta_{i,n_i,m}\varepsilon]\}\\
=&\{\bold{z}_i+\vec{\delta}_{i,1}\varepsilon,\ldots,\bold{z}_i+\vec{\delta}_{i,n_i}\varepsilon\},
\end{aligned}
\end{scriptsize}
\end{equation}
where  $|\delta_{i,j,r}|<1$  for all $i=1,\ldots,k$,
$j=1,\ldots,n_i$, $r=1,\ldots,m$. Let $U_k$ be the partially row
reduced form obtained by   applying $k$ steps of the GECP to  the
matrix of traces $R$ corresponding to $C_1\cup \cdots\cup C_k$. Then
$R$ and $U_k$ have entries from the field $\mathbb{C}(\varepsilon)$.
\end{defn}

\begin{assu}\label{perm}
 Based on Proposition \ref{regular},  we will assume  that if
the GECP applied to $R$ produces the permutation matrices $P$ and
$Q$ then the matrix $P\,R_0\,Q$ is regular, where
$R_0=R|_{\varepsilon=0}$. To simplify the notation for the rest of
the paper we will assume that $Q=id,$ i.e. the rows and columns of
$PRQ=PR$ correspond to the bases
\begin{eqnarray}\label{sigmaB}
\sigma B=[b_{\sigma(1)}, \ldots, b_{\sigma(n)}]\;\text{ and }
\;B=[b_1, \ldots,b_n],
\end{eqnarray}
respectively, where $\sigma$ is the permutation corresponding to
the matrix $P$. This assumption does not constrain the generality
since we may rename $B$ in the definition of $R$.
\end{assu}

With the assumption that $P\,R_0$ has rank $k$ and is regular, we
can assume that all the denominators appearing in the entries of
$U_k$ are minors of $R$ which
 are non-zero at $\varepsilon=0$. Therefore we
can take their Taylor expansion around $\varepsilon=0$ and
consider them as elements  of the formal power series ring
$\mathbb{C}[[\varepsilon]]$. In this ring we shall work with
residue classes modulo $\varepsilon^2$, i.e., in some
considerations we factor out the ideal
$\langle\varepsilon^2\rangle$ of
$\mathbb{C}[[\bold{\varepsilon}]]$.

The results in the rest of the paper are all valid modulo
$\varepsilon^2$ in the formal power series setting described above.
In practice what this means is that the method we propose works up
to a precision which is the double of the original size of the
clusters.

\begin{rem}\label{notall}\emph{
In Definition \ref{setting} we assume that the clusters are linear
perturbations of a set of multiple roots.  Note that not all
multiplicity structures can be obtained as a limit of such clusters
with linear perturbation of fixed directions $\vec{\delta}_{i,j}$
 as $\varepsilon$ approaches $0$. However, as we have seen in
 Proposition \ref{traceroots},  the matrix of traces at
 $\varepsilon=0$ is independent of the directions
 $\vec{\delta}_{i,j}$, and in fact does not depend on the multiplicity
 structure of the roots. Since all the subsequent results in the paper only depend on the matrix of traces
 and are
 only valid modulo $\varepsilon^2$, we do not limit the generality
 by considering only linear perturbations.  This is not true however for the
 multiplication matrices, which depend on the
 multiplicity structure  of the roots at $\varepsilon=0$, as seen in Example \ref{clusterex}.}
\end{rem}

In order to describe the structure of the matrices in the $LU$
decomposition of the matrix of traces obtained by GECP in terms of
the elements in the clusters, we need the following definition:

\begin{defn} \label{vander}
Let $B=[b_1,\ldots,b_n]\in \mathbb{C}[x_1, \ldots, x_m]^n$, and
let $\bold{z}_1, \ldots,  \bold{z}_r\in  \mathbb{C}^m$ be not
necessary distinct points. We call the $n\times r$ matrix
$$
V:= \left[ b_i(\bold{z}_j)\right]_{i,j=1}^{n,r}
$$
the {\it Vandermonde matrix} of $\bold{z}_1, \ldots, \bold{z}_r$
w.r.t. $B$. Note that if $r=n$ then the matrix of traces in
Definition \ref{rmat} and the Vandermonde matrix are closely
related:
$$
R= VV^T.
$$
\end{defn}

 The following proposition
gives asymptotic bounds for the entries of the matrix obtained from
a partial Gauss elimination with complete pivoting on the matrix of
traces $R$ for the case where the $n$ roots of the system correspond
to $k$ clusters, each of them with $n_i$ roots ($i=1,\ldots,k$) and
radius proportional to $\varepsilon$ in the max-norm.

\begin{prop}\label{elemrow}
Let $B=[b_1,\ldots,b_n]\in \mathbb{C}[x_1, \ldots, x_m]^n$. Let
$\{\bold{z}_1, \ldots, \bold{z}_k\}\in\mathbb{C}^m$ and  the
clusters
 $C_1, \ldots,
C_k$ around $\{\bold{z}_1, \ldots, \bold{z}_k\}$ be as in
Definition \ref{setting}.

Let $R$ be the matrix of traces associated to $C_1\cup\cdots\cup
C_k$ and $B$ (see   Definition \ref{tracemat} and Proposition
\ref{traceroots}). Let  $P$ and $R_0:=R|_{\varepsilon=0}$ be as in
 Assumption \ref{perm} and assume that $P\,R_0$
 has rank $k$ and is regular. Then, after $k$
steps of the regular Gaussian elimination on $P\,R$ we get a
partially row reduced matrix $U_k$, such that its last $n-k$ rows
satisfy
\begin{small}
\begin{eqnarray}\label{Taylor}
[ U_k]_{i,j}=
\begin{cases}
 0, &\text{if }j\leq k\\
c_{i,j}\varepsilon^2 + h. o. t.(\varepsilon)\;\;\in
\mathbb{C}[[\varepsilon]]&\text{if }j>k
\end{cases} \;\;\text{ for }\;\;i=k+1, \ldots, n.
\end{eqnarray}
\end{small}
The values of $c_{i,j}\in \mathbb{C}$ depends on $\ n$,
$\{\bold{z}_1, \ldots,\bold{z}_k\}$, $\{\vec{\delta}_{s,t}\}$ and
$B$ (we will give a bound for $c_{i,j}$ in  Proposition
\ref{bound}). Here $h.o.t.(\varepsilon)$ denotes the higher order
terms in $\varepsilon$. Moreover, the formal power series in
(\ref{Taylor}) are convergent in a sufficiently small neighborhood
of $\varepsilon=0$.
\end{prop}

\begin{proof}
To simplify the notation, denote $\tilde{R}=P\,R$. The proof is
based on the fact that after $k$ steps of the regular Gaussian
elimination on $\tilde R$, the partially reduced $U_k$ matrix has
elements $(i,j)$, for $i,j=k+1,\ldots,n$, of the form
\begin{eqnarray}\label{Uij}
\frac{\det \,(\tilde{R}^{(k+1)}_{i,j})}{\det\,( \tilde{R}^{(k)})}
\end{eqnarray}
 where $\tilde{R}^{(k)}$ is the $k\times k$ principal submatrix of $\tilde{R}$
 and $\tilde R^{(k+1)}_{i,j}$ is the  $(k+1)\times (k+1)$ submatrix
of $\tilde R$ corresponding to rows $\{1,\ldots,k,i\}$ and columns
$\{1,\ldots,k,j\}$. This follows at once from the facts that both
the numerator and the denominator of (8) stay the same during the
row operations performed, and the reduced form of
$R_{i,j}^{(k+1)}$ is upper triangular.

Let $V$ be the $n\times n$ Vandermonde matrix of
$C_1\cup\cdots\cup C_k$ with respect to $B$ and recall that $R=V\,
V^T$, thus $\tilde R = (PV)(V^T)$. Let $\sigma$  be the
permutation corresponding to $P$ and let $\sigma B$ be as in
(\ref{sigmaB}). Observe that  $$\tilde R^{(k+1)}_{i,j}=V_{\sigma
B_i}\,V_{ B_j}^T,$$ where $V_{\sigma B_i}$ and $V_{B_j}$  are the
$(k+1) \times n$ Vandermonde matrices corresponding to
$C_1\cup\cdots\cup C_k$ and respectively
  to $\sigma B_i:=[b_{\sigma(1)}, \ldots,
b_{\sigma(k)}, b_{\sigma(i)}]$,
 and $B_j:=[b_1, \ldots, b_k, b_j]$.
 Therefore, by the Cauchy-Binet formula we have
 \begin{eqnarray}\label{subdet}
 \det ( \tilde R^{(k+1)}_{i,j}) =\sum_{|I|=k+1} \det (V_{\sigma B_i, I})\det(V_{B_j, I}),
 \end{eqnarray}
where $V_{\sigma B_i, I}$ denotes the $(k+1)\times (k+1)$
submatrix of $V_{\sigma B_i}$ with columns corresponding to the
points in $I\subset C_1\cup\cdots\cup C_k$, and the summation is
taken for all $I\subset C_1\cup\cdots\cup C_k$ such that
$|I|=k+1$. Note that all the determinants in (\ref{subdet}) are
polynomials in $\varepsilon$. Since ${\rm rank}\,(
V|_{\varepsilon=0})=k$, we have $\det\, (V_{\sigma B_i, I}
)|_{\varepsilon=0}=\det\, (V_{ B_j, I} )|_{\varepsilon=0}=0$, thus
they are divisible by $\varepsilon$ for all $i=k+1, \dots n$ and
$I \subset C_1\cup\cdots\cup C_k $ with $|I|=k+1$.  Therefore we
get that $\det \,( \tilde R^{(k+1)}_{i,j})$ is divisible by
$\varepsilon^2$.

Finally we note that the assumption that $P\, R_0$ has rank $k$ and
is regular implies that $$\det(\tilde R^{(k)}|_{\varepsilon=0})\neq
0,$$ which proves that the Taylor expansion of the ratio in
(\ref{Uij}) around $\varepsilon=0$ has zero constant and linear
terms, as was claimed. The formal power series in (\ref{Taylor}) are
convergent in a sufficiently small neighborhood of $\varepsilon=0$,
since they are the Taylor series of rational functions with non-zero
denominators at $\varepsilon=0$.
\end{proof}

From the previous results it follows that if we have $k$ clusters
of size $n_i$, with $i=1,\ldots,k$, $\sum_{i=1}^kn_i=n$, then
after $k$ steps of GECP on the matrix of traces $R$, we get the
matrix

\begin{small}
\begin{eqnarray}\label{LUshape}
U_k= \left[
\begin{array}{cccccc}
[U_k]_{1,1}&&\cdots&\cdots&\cdots&
[U_k]_{1,n}\\
0&\ddots&\cdots&\cdots&\cdots&\vdots\\
&&[U_k]_{k,k}&\cdots&\cdots&[U_k]_{k,n}\\
\vdots&&0 &c_{k+1,k+1}\varepsilon^2&\cdots&c_{k+1,n}\varepsilon^2\\
&&\vdots&\vdots&\ddots&\vdots\\
0&&0 &c_{n,k+1}\varepsilon^2&\cdots&c_{n,n}\varepsilon^2\\
\end{array}
\right ]+ h.o.t.(\varepsilon)
\end{eqnarray}
\end{small}
where the constant term in $\varepsilon$ of $[U_k]_{i,i}$ is
non-zero for $i\leq k$.

The next proposition gives a bound for the coefficient $c_{i,j}$
of $\varepsilon^2$ in (\ref{LUshape}). It also gives an idea of
the magnitude of the threshold one can use to decide on the
numerical  rank which would additionally indicate how small the
size of the clusters need to be for our method to work.

\begin{prop}\label{bound}
Let $B=[b_1,\ldots,b_n]\in \mathbb{C}[x_1, \ldots, x_m]^n$. Let
$\{\bold{z}_1, \ldots, \bold{z}_k\}\in\mathbb{C}^m$. Let the
clusters
 $C_1, \ldots,
C_k$ around $\{\bold{z}_1, \ldots, \bold{z}_k\}$ be as in
(\ref{cluster}) with $|\delta_{i,j,r}|\leq 1$ for all
$i=1,\ldots,k$, $j=1,\ldots,n_i$, $r=1,\ldots,m$. Let $R$ be the
matrix of traces associated to $C_1\cup\cdots\cup C_k$ and $B$.

Let $\;\bold{b'}$ be such that
\begin{eqnarray}\label{b'}
\bold{b'}\geq max_{\{l,i,r\}}\left\{\left|\frac{\partial
b_{l}}{\partial x_{r}}(\bold{z}_i)\right|\right\}. \end{eqnarray}

Assume that the GECP applied to $R$ also implies complete pivoting
on $\left.R\right|_{\varepsilon=0}$. Then the bound for the
coefficients $c_{i,j}$ of $\varepsilon^2$ in the $U_k$ matrix,
obtained after $k$ steps of the  GECP applied to the matrix of
traces $R$, is given by
\[
|c_{i,j}|\leq \alpha\cdot (\bold{b'})^2,\] where
$\alpha=4(n-k)(k+1)^2m^2$.
\end{prop}


\begin{proof}
We denote  $\tilde{R}:=P\,R$, where $P$ is a permutation matrix
such that the first $k$ steps of GECP applied to both $P\,R$ and
$P\,R|_{\varepsilon=0}$ is well defined and the same as regular
Gaussian elimination. Note that we need the assumption that GECP
applied to $R$ also implies complete pivoting on
$\left.R\right|_{\varepsilon=0}$ since Proposition \ref{regular}
only implies that $P\,R|_{\varepsilon=0}$ is regular, but below we
will  also need the pivots in $P\,R|_{\varepsilon=0}$ to have
maximal absolute values. One can achieve this by making the right
selection among equal possible pivots while performing GECP on
$R$. We will use this assumption at the end of the proof.

Denote the bases corresponding to the rows and columns of $\tilde R$
by $\sigma B=[b_{\sigma(1)}, \ldots, b_{\sigma(n)}]\;\text{ and }
\;B=[b_1, \ldots,b_n]$ as in (\ref{sigmaB}).

The partially reduced $U_k$ matrix has elements $(i,j)$, for
$i,j=k+1,\ldots,n$, of the form
\begin{eqnarray}\label{Uij}
\frac{\det \,(\tilde{R}^{(k+1)}_{B',B''})}{\det\,( \tilde{R}^{(k)})}
\end{eqnarray}
where $\tilde{R}^{(k)}$ is the $k$-th principal submatrix of $\tilde
R$ and $\tilde R^{(k+1)}_{B',B''}$ is the  $(k+1)\times (k+1)$
submatrix of $\tilde R$ corresponding to rows
$B':=[b_{\sigma(1)},\ldots, b_{\sigma(k)},b_{\sigma(i)}]$ and
columns  $B'':=[b_1, \ldots, b_k, b_j]$. In order to get an upper
bound for $|c_{i,j}|$, we will get an upper bound for the
coefficient
 of $\varepsilon^2$ in $\det\,(\tilde R^{(k+1)}_{B ',B''})$ and divide it by  the constant term
 of $|\det\,(\tilde R^{(k)})|$.

 Fix $i, j\in \{1, \ldots, n\}$. We will use the Cauchy-Binet
 formula
\[
 \det ( \tilde R^{(k+1)}_{B',B''}) =\sum_{|I|=k+1} \det (V_{ B', I})\det(V_{B'',
 I}),
\]
 where the summation is for
 $I\subset C_1\cup\cdots\cup C_k$  of cardinality $k+1$, and $V_{ B', I}$ and $V_{B'',
 I}$ are the Vandermonde matrices corresponding to $I$   w.r.t. $B'$ and $B''$, respectively.
 Since the derivative of the determinant of a matrix is the sum of determinants
obtained by
 replacing one by one the columns  by the derivative of that column, after expanding the
 determinants in the sum
 by
 their columns containing the derivatives, we get

 \[
\left.\frac{\partial \det \,( V_{ B',
I})}{\partial\varepsilon}\right|_{\varepsilon=0} =\sum_{b'\in
B'}\sum_{{\bf z}\in I}\pm\left(\sum_{t=1}^m\delta_{{\bf
z},t}\frac{\partial b'}{\partial x_t}({\bf
z})|_{\varepsilon=0}\right)\det \,( \left.V_{B'-\{b'\}, I-\{{\bf
z}\}})\right|_{\varepsilon=0},
\]
where $\delta_{{\bf z},t}$ is the coefficient of $\varepsilon$ in
the $t$-th coordinate of ${\bf z}\in I$. We can obtain a similar
expression for $\left.\frac{\partial \det \,( V_{B'',
I})}{\partial\varepsilon}\right|_{\varepsilon=0}$.

 Note that $\det \,( \left.V_{B'-\{b\}, I-\{{\bf
z}\}})\right|_{\varepsilon=0}$ is non-zero only if
$I|_{\varepsilon=0}=\{{\bf z}_1, \ldots, {\bf z}_k\}\cup\{{\bf
z}_i\}$ for some $i=1,\ldots,k$ and ${\bf z}|_{\varepsilon=0}={\bf
z}_i$. In that case $I-\{{\bf z}\}|_{\varepsilon=0}=\{{\bold
z}_1,\ldots,{\bold z}_k\}$, which we denote by
$$\zz:=\{{\bold z}_1,\ldots,{\bold z}_k\}$$ for simplicity.

Thus we have that if $I|_{\varepsilon=0}=\{{\bf z}_1, \ldots, {\bf
z}_k\}\cup\{{\bf z}_i\}$ then

 \[
 \begin{aligned}
\left|\frac{\partial \det \,( V_{ B',
I})}{\partial\varepsilon}\right|_{\varepsilon=0}\leq\sum_{b'\in
B'}\sum_{{\bf z}\in I\atop{{\bf z}|_{\varepsilon=0}={\bf z}_i}}
\left|\det \,( V_{B'-\{b'\}, \zz})\right|
\sum_{t=1}^m\left|\frac{\partial b'}{\partial x_t}({\bf
z})\right|_{\varepsilon=0}
\end{aligned}
\]
using that $|\delta_{{\bf z},t}|\leq 1$. Therefore, we get that the
coefficient of $\varepsilon^2$ in $\det \,(\tilde
R^{(k+1)}_{B',B''})$ is bounded by
  \[
  \begin{aligned}
4\sum_{b'\in B'\atop{b''\in B''}}&\left|\det \,( V_{B'-\{b'\},
\zz})\det \,(
V_{B''-\{b''\}, \zz})\right|\notag\\
&\left(\sum_{|I|=k+1\atop{I|_{\varepsilon=0}=\{{\bold
z}_1,\ldots,{\bold z}_k\}\cup\{{\bold
z}_i\}}}\left(\sum_{t=1}^m\left|\frac{\partial b'}{\partial
x_t}({\bf z}_i)\right|\right)\left(\sum_{t=1}^m\left|\frac{\partial
b''}{\partial x_t}({\bf z}_i)\right|\right)
\right).\notag\\
\end{aligned}
\]
using the fact that there are two possible ways to pick ${\bf z}\in
I$ with ${\bf z}|_{\varepsilon=0}={\bf z}_i$ from
$I|_{\varepsilon=0}=\{{\bold z}_1,\ldots,{\bold z}_k\}\cup\{{\bold
z}_i\}$.

Using the upper bound $\bold{b'}$ and counting the number of times
we can choose $I\subset C_1\cup\cdots\cup C_k$ such that $|I|=k+1$
and $I|_{\varepsilon=0}=\{{\bold z}_1,\ldots,{\bold
z}_k\}\cup\{{\bold z}_i\}$, we get
  \[
4m^2\bold{b'}^2(n-k)\left(\prod_{j=1}^k n_j\right)\left(\sum_{b'\in
B'\atop{b''\in B''}}\left|\det \,( V_{B'-\{b'\}, \zz})\det \,(
V_{B''-\{b''\}, \zz})\right|\right).
\]

On the other hand, we have that if $ \tilde R_{B'-\{b'\},
B''-\{b''\}}$ is the matrix of traces with rows corresponding to the
$B'-\{b'\}$ and columns corresponding to $B''-\{b''\}$ then
\[
\left.\det \,(\tilde R_{B'-\{b'\},
B''-\{b''\}})\right|_{\varepsilon=0}=\left(\prod_{i=1}^k
n_i\right)\det \,( V_{B'-\{b'\}, \zz})\det \,( V_{B''-\{b''\},
\zz}).
\]

Therefore, the bound for the coefficient of $\varepsilon^2$ in $\det
\,(\tilde R^{(k+1)}_{B',B''})$ is

\begin{eqnarray}\label{partial2}
4m^2\bold{b'}^2(n-k)\sum_{b'\in B'\atop{b''\in B''}}\left|\det
\,(\tilde R_{B'-\{b'\}, B''-\{b''\}})\right|_{\varepsilon=0}.
 \end{eqnarray}

Next we use the assumption above on $\tilde
R|_{\varepsilon=0}=PR|_{\varepsilon=0}$ to have maximal pivots in
the first $k$ diagonal entries to get
 \begin{eqnarray}\label{kk}
 |\det\,(\tilde R^{(k)})|_{\varepsilon=0}\geq
|\det \,(\tilde R_{B'-\{b'\}, B''-\{b''\}})|_{\varepsilon=0}
 \end{eqnarray}
 which is true since  the left and right hand side of (\ref{kk}) divided by
$|\det\,(\tilde R^{(k-1)})|_{\varepsilon=0}$ give the absolute
values of the entries of the partially row reduced matrix after
$k-1$ steps of GECP. Therefore we can replace $|\det \,(\tilde
R_{B'-\{b'\}, B''-\{b''\}})|_{\varepsilon=0}$ by $ |\det\,(\tilde
R^{(k)})|_{\varepsilon=0}$ in (\ref{partial2}) and divide the
expression (\ref{partial2}) by $|\det\,(\tilde
R^{(k)})|_{\varepsilon=0}$, thus  we get the following bound for
the coefficient $c_{i,j}$ of $\varepsilon^2$ in the $U_k$ matrix:
\[
4m^2\bold{b'}^2(n-k)(k+1)^2.
\]
\end{proof}

\begin{rem}\emph{ The above proposition gives estimates in terms of $\{\bold{z}_1, \ldots, \bold{z}_k\}$,
which we do not assume to know a priori. The following heuristic
methods can be used to check whether the estimated numerical rank
is correct, given a required precision $\varepsilon$. Assuming
that we know the magnitude of the coordinates of the roots, we may
compute the matrix of traces corresponding to $n$ randomly chosen
distinct roots which have the same order of magnitude as the
original roots. Then comparing the diagonal entries of the
$U$-matrices obtained by applying the GECP for the matrices of
traces, we can set the numerical rank to be the first entry where
the discrepancy is of order $\varepsilon^2$. Another heuristics is
to increase $k$ one by one, compute the approximate radical ideal
(see Definition \ref{apprad}) corresponding to the case when the
numerical rank of $R$ is $k$. Compute the roots of the approximate
radical ideal, and substitute them back into the original system.
If the error is of order of magnitude $\varepsilon$, accept $k$
and the computed approximate radical as the output.}
\end{rem}

\begin{ex} \label{ex3cont} Continuing Example \ref{ex3}, we apply
the GECP to the matrix $R$ defined in (\ref{rmatc}). After two
steps of GECP we obtain the following matrix: {\scriptsize$$
U_2=\left[
\begin {array}{ccccc} 11.07093& -5.04240& 7.03226& -1.01900&
-1.04951\\0& 8.71265& 2.18381& 6.53716& 6.55387\\0& 0&
0.454213\times 10^{-4}& 0.7407\times 10^{-5}& 0.178036\times
10^{-3}\\0& 0& 0.7397\times 10^{-5}& 0.728\times 10^{-6}&
0.41955\times 10^{-4}\\0& 0& 0.188071\times 10^{-3}& 0.52002\times
10^{-4}& 0.657084\times 10^{-3}
\end{array} \right]
$$}
with columns permuted so that they correspond to the basis
$[x_1x_2, x_2, x_1, 1, x_1^2]$. Note that the largest entry in the
$3\times 3$ bottom right corner of $U_2$ is  between $\varepsilon$
and $\varepsilon^2$ (here $\varepsilon\approx 0.01$ in this
example). Thus we consider the numerical rank of $R$ to be $2$.
From the nullspace of the first two rows of $U_2$ we can obtain
the following approximate multiplication matrices:
 {\scriptsize$$  \mathcal{M}'_{x_1} =
    \left[ \begin {array}{cc}
 0 & 1\\
1.00382& -0.37849\times 10^{-3}
\end {array} \right]  \quad \mathcal{M}'_{x_2} =
    \left[ \begin {array}{cc}
 1.49973& -0.49972\\
-0.5016325& 1.50162
\end {array} \right]
$$}
(see Section \ref{AppRad} below for more details on approximate
multiplication matrices). The eigenvalues of $\mathcal{M}'_{x_1}$
and $\mathcal{M}'_{x_2}$ are respectively {\small
 $$ \{1.000018,
-1.003803\}\;\;\text{  and }\;\;\{ 0.9999943, 2.001349\}.
 $$}
 Note that these eigenvalues are close to the avarages of the
 coordinates of the roots in  the two clusters.
\end{ex}

\section{Singular Values of $R$}

Using the previous results we will now study  the singular values
of the matrix of traces $R$ of a system with clusters of roots. We
denote  $\tilde{R}:=P\,R$, where $P$ is a permutation matrix
obtained by $k$ steps of GECP applied to $R$ and we assume that
$P\,R|_{\varepsilon=0}$ is  regular,  as in Assumption \ref{perm}.
Let $U_k$ be the matrix obtained after $k$ steps of GECP on the
matrix $\tilde{R}$, as in (\ref{LUshape}). Let $\widehat{U}_k$ be
the matrix obtained after replacing the last $n-k$ rows of $U_k$
by zeros. Let $L_k$ be such that $\tilde{R}=L_kU_k$ (in other
words $L_k$ is the transformation matrix obtained after $k$ steps
of GECP on $\tilde{R}$). Let $\widehat{R}=L_k\widehat{U}_k$. Using
the submultiplicative property of matrix norms, we have that
\[
\|\tilde{R}-\widehat{R}\|_F=\|L_kU_k-L_k\widehat{U}_k\|_F\leq
\|L_k\|_F\|U_k-\widehat{U}_k\|_F
\]
where $\|\cdot\|_F$ denotes the Frobenius matrix norm.

Let $\sigma_1\geq\cdots\geq\sigma_n$ be the singular values of
$R$, which are also the singular values of $\tilde{R}$. Since by
definition $\sigma_i$ is the 2-norm distance from $\tilde{R}$ to
the nearest rank $i$ matrix, and $\widehat{R}$ is an $n$ by $n$
matrix of rank $k$, we have that
\[
\sigma_n\leq\cdots\sigma_{k+1}\leq\|\tilde{R}-\widehat{R}\|_2.
\]

Given that the 2-norm of a matrix is smaller than or equal to its
Frobenius norm, we have
\[
\sigma_n\leq\cdots\sigma_{k+1}\leq
\|L_k\|_F\|U_k-\widehat{U}_k\|_F.
\]

Since we are using GECP it is easy to see that
\[
\left[L_k\right]_{i,j}\leq 1\,\text{ for all } j=1,\ldots,k,\, i>j
\]

and the matrix $L_k$ obtained after $k$ steps of GECP is of the
form
\begin{small}
\[
L_k= \left[
\begin{array}{ccccccc}
1&0&\cdots&\cdots&\cdots&\cdots&0\\
\left[L_k\right]_{2,1}&\ddots & \ddots &\cdots & \cdots & \cdots & \vdots\\
\vdots&\ddots&1&0&\cdots&\cdots&\vdots\\
\vdots&\cdots&\left[L_k\right]_{k+1,k}&1&0&\cdots&\vdots\\
\vdots&\cdots&\vdots&0&\ddots&\ddots& \vdots\\
\vdots&\cdots&\vdots&\vdots&\ddots&\ddots&0\\
\left[L_k\right]_{n,1}&\cdots&\left[L_k\right]_{n,k}&0 &\cdots&0&1\\
\end{array}
\right ].
\]
\end{small}

Therefore we have
\[
\|L_k\|_F\leq\sqrt{\frac{2n+2nk-k^2-k}{2}}.
\]

From Proposition \ref{bound} we have that for $i,j=k+1\ldots n$,
the elements of  $U_k$ are of the form
\[
[U_k]_{ij}=\omega\varepsilon^2+h.o.t.(\varepsilon),
\]
where $\omega=4(n-k)(k+1)^2m^2 \bold{b'}^2$ and $\bold{b'}^2$ is
defined in (\ref{b'}).

We therefore have
\[
\|U_k-\widehat{U}_k\|_F=\sqrt{\sum_{i,j=k+1}^n([U_k]_{ij})^2}\leq
(n-k)\omega\varepsilon^2+h.o.t.(\varepsilon).
\]

We summarize the above argument in the next Proposition, showing
that the $k+1$-th singular value of $R$ is asymptotically equal to
$\varepsilon^2$.

\begin{prop}\label{SVD}
 Let $R$ be the
matrix of traces associated to $C_1\cup\cdots\cup C_k$ and $B$
where the clusters $C_1, \ldots, C_k$ around $\bold{z}_1, \ldots,
\bold{z}_k\in\mathbb{C}^m$ are
 as in (\ref{cluster}) and  $B=[b_1,\ldots,b_n]\in \mathbb{C}[x_1, \ldots, x_m]^n$.
  Let $\sigma_1\geq\cdots\geq\sigma_n$ be
the singular values of $R$. Then
$$
\sigma_{k+1}= \Omega\varepsilon^2 + h.o.t.(\varepsilon)
$$
where
$$
\Omega \leq
4(n-k)^2(k+1)^2m^2\sqrt{\frac{2n+2nk-k^2-k}{2}}(\bold{b'})^2
$$
and $\bold{b'}$ is defined in (\ref{b'}).
\end{prop}

\begin{ex}
Continuing Example \ref{ex3}, we compute the singular values of
the matrix $R$ defined in (\ref{rmatc}):
\[[22.8837, 14.2433, 0.448334\times 10^{-3}, 0.174904\times 10^{-4}, 0.594796\times
10^{-5}].
\]
 We  have that the third singular value is between $\varepsilon$
 and
$\varepsilon^2$ (in this example $\varepsilon\approx 0.01$), thus
we can set  the numerical rank of the matrix $R$ to be 2. Note
that the 2-norm distance of the matrix $R$ from
$R|_{\varepsilon=0}$ is not the same order of magnitude as the
third singular value, it is $0.147$ as was computed in Example
\ref{ex3}. This is the reason why
 we used the partial LU-decomposition of $R$  and not $R|_{\varepsilon=0}$ to obtain a bound
 for $\sigma_{k+1}$.
\end{ex}

\section{Approximate Radical Ideal}\label{AppRad}

Using our previous results, we can now define the concept of an
approximate radical ideal and describe its roots in terms of the
elements of the clusters.

\begin{defn}\label{apprad}
Let $B=[b_1,\ldots,b_n]\in \mathbb{C}[x_1, \ldots, x_m]^n$ and the
clusters
 $C_1, \ldots,
C_k$  be as in Definition \ref{setting}.
 Let
$R$ be the matrix of traces associated to $C_1\cup\cdots\cup C_k$
and $B$.  Let the permutation matrix $P$  corresponding to the
permutation $\sigma$ obtained after $k$ steps of GECP on $R$ as in
 Assumption \ref{perm}, so that the rows and columns of
$\tilde{R}:=P\,R$ correspond to $\sigma B$ and $B$, respectively,
as in (\ref{sigmaB}).
 We define the vectors ${\bf v}_{i,j}\in \mathbb{C}(\varepsilon)^{k}$
 for $ i=1,
\ldots m $ and $ j=1, \ldots k$, as the solutions of the following
$mk$ linear systems:
 \begin{equation}\label{defeq}
\tilde{R}^{(k)} {\bf v}_{i,j}  = {\bf r}_{i,j} \quad i=1, \ldots,
m \quad j=1, \ldots, k,
 \end{equation}
 where the left hand sides are always  the $k\times k$ principal submatrix of $\tilde{R}$,
while for any fixed $i$ and $j$ the right hand side of
(\ref{defeq}) is defined as
\begin{equation}\label{righthandside}
{\bf r}_{i,j}:= \left[
\begin{array}{c}Tr(x_ib_jb_{\sigma(1))}\\\vdots\\Tr(x_ib_jb_{\sigma(k)})
\end{array}\right]\;\;\in \;\;\mathbb{C}^k.
\end{equation}
Note that one can compute the vectors ${\bf r}_{i,j}$ the same way
as the columns of the matrix of traces. Then we define the
following $mk$ polynomials:

\begin{eqnarray}\label{fij} f_{i,j}:=x_ib_{j}-\left(\sum_{s=1}^k [{\bf v}_{i,j}]_{s} b_{s} \right)\quad
i=1, \ldots, m, \quad j=1, \ldots, k. \end{eqnarray}

  We will call the {\it approximate radical ideal} of the clusters $C_1\cup \cdots \cup
C_k$
 the ideal generated by
$$\widetilde{\sqrt{I}}:=\langle f_{i,j}\;:\;i=1, \ldots, m, \; j=1, \ldots, k\rangle.
 $$
  We also
define the {\em approximate multiplication matrices of the
radical} of $C_1\cup\cdots\cup C_k$ with respect to the basis
$[b_{1}, \ldots, b_{k}]$ to be the matrices ${M}'_{x_1},
\ldots,{M}'_{x_m}\in\mathbb{C}(\varepsilon)^{k\times k}$ where
$$
[M'_{x_i}]_{j,s}:= [{\bf v}_{i,j}]_{s}\quad \;i=1, \ldots, m, \;
j,s=1, \ldots, k.
$$
\end{defn}

\begin{rem}\emph{
We can also define the approximate multiplication matrices of the
radical of $C_1\cup\cdots\cup C_k$ from a system of multiplication
matrices of $C_1\cup\cdots\cup C_k$ with respect to $B$ by
changing the basis as follows: Let $r_{k+1}, \ldots, r_n\in
\mathbb{C}(\varepsilon)^n$ be a basis for the nullspace of the
first $k$ rows of $P\,R $. Let $v_1, \ldots, v_k\in \mathbb{C}^n$
be such that  $B':=[v_1, \ldots, v_k, r_{k+1}, \ldots, r_n]$ forms
a basis for $\mathbb{C}(\varepsilon)^n$. Let ${M}_{x_1},
\ldots,{M}_{x_m}\in \mathbb{C}(\varepsilon)^{n\times n}$ be the
multiplication matrices of the clusters $C_1\cup\cdots\cup C_k$
with respect to the basis $B'$. Then  the
 approximate multiplication matrices of the radical of
$C_1\cup\cdots\cup C_k$ with respect to $[v_1, \ldots, v_k]$ are
the matrices ${M}'_{x_1},
\ldots,{M}'_{x_m}\in\mathbb{C}(\varepsilon)^{k\times k}$ obtained
as the principal $k\times k$ submatrices of ${M}_{x_1},
\ldots,{M}_{x_m}$, respectively. Note that the eigenvalues of
${M}_{x_i}$ are the $x_i$ coordinates of the elements of the
clusters reordered in a way that the first $k$ correspond to one
eigenvalue from each cluster. However, we also remark that  we
have to be careful with the multiplication matrices ${M}_{x_1},
\ldots,{M}_{x_m}$ since they are not always continuous at
$\varepsilon=0$, as noted in Remark \ref{notall}, thus we cannot
consider their entries as elements of $\mathbb{C}[[\varepsilon]]$.
That is the reason we chose to define the approximate radical as
in Definition \ref{apprad}.}
\end{rem}

The next proposition asserts that  when $\varepsilon=0$ our
definition gives the multiplication matrices of the radical ideal.

\begin{prop}\label{continuous}
Using the assumptions of Definition \ref{apprad}, the coordinates
of the vectors ${\bf v}_{i,j}\in \mathbb{C}(\varepsilon)^{k}$
defined in (\ref{defeq}) are continuous in $\varepsilon=0$ for all
$i=1, \ldots, m$ and  $j=1, \ldots, k$. Furthermore, the points
${\bf z}_1, \ldots, {\bf z}_k$ are common roots of the polynomials
 $\{ f_{i,j}|_{\varepsilon=0}\;:\;i=1, \ldots, m, \;
j=1, \ldots, k\}$, and
   the matrices
$${M}'_{x_1}|_{\varepsilon=0}, \ldots, {M}'_{x_m}|_{\varepsilon=0}$$
form a system of multiplication matrices for the algebra
$\mathbb{C}[{\bf x}]/\sqrt{I} $.
\end{prop}

\begin{proof}
Using Assumption \ref{perm}, the continuity of the coordinates of
the vectors ${\bf v}_{i,j}\in \mathbb{C}(\varepsilon)^{k}$ follows
from our assumption that the $k\times k$ principal submatrix
$\tilde{R}^{(k)}$ of $\tilde{R}$ is nonsingular at
$\varepsilon=0$.

Next we  show that ${\bf z}_1, \ldots, {\bf z}_k$ are roots of
$f_{i,j}|_{\varepsilon=0}$ for all $i\in \{1, \ldots, m\}$ and
$j\in \{1, \ldots, k\}$. Fix $i$ and $j$. Assume that
\begin{eqnarray}\label{inter}
x_ib_{j}-\left(\sum_{s=1}^k w_{i,j,s} b_s\right)= 0
 \end{eqnarray}
  is satisfied by ${\bf z}_1, \ldots, {\bf z}_k$, which is equivalent to the
column vectors
 \begin{eqnarray}\label{cij} {\bf w}_{i,j}:=\left[w_{i,j,1}, \ldots,
w_{i,j,k}, -1\right]^T
\end{eqnarray}
satisfying the homogeneous linear system with coefficient matrix
$W$ defined to be the transpose of the $(k+1)\times k$ Vandermonde
matrix of ${\bf z}_1, \ldots, {\bf z}_k$ with respect to $[b_1,
\ldots, b_k,x_ib_{j}]$.

On the other hand, by (\ref{defeq}), the vector $[{\bf
v}_{i,j}|-1]_{\varepsilon=0}$ is in the nullspace of the $k\times
(k+1)$ matrix $[\tilde{R}^{(k)}|{\bf r}_{ij}]_{\varepsilon=0}$. We
have
$$[\tilde{R}^{(k)}|{\bf r}_{ij}]_{\varepsilon=0}=V_1V_2^T$$
 where $V_1$ and $V_2$ are the Vandermonde
matrices of $C_1, \ldots, C_k$ at $\varepsilon=0$ corresponding
respectively  to $[b_{\sigma(1)}, \ldots, b_{\sigma(k)}]$ and
$[b_1, \ldots, b_k,x_ib_{j}]$, thus $V_2^T$ is the same as $W$
except the row corresponding to ${\bf z}_s$ is repeated $n_s$
times for $s=1, \ldots, k$. This implies that  the nullspace of
$W$ is a subset of the nullspace of $[\tilde{R}^{(k)}|{\bf
r}_{ij}]_{\varepsilon=0}$. But since both nullspaces has dimension
one, we must have ${\bf w}_{i,j}=[{\bf
v}_{i,j}|-1]_{\varepsilon=0}$, i.e. $f_{i,j}|_{\varepsilon=0}=0$
is satisfied by ${\bf z}_1, \ldots, {\bf z}_k$.

 Next we prove that the matrices
$M'_{x_1}|_{\varepsilon=0}, \ldots, M'_{x_d}|_{\varepsilon=0}$
form a system of multiplication matrices for $\mathbb{C}[{\bf
x}]/\sqrt{I}$. First note that for any $g\in \mathbb{C}[{\bf x}]$,
if ${\bf z}$ is a common root of the system
$$
gb_j-\sum_{s=1}^k c_{j,s} b_s=0 \quad j=1, \ldots, k
$$
and ${\bf z}$ is not a common root of $b_1, \ldots, b_k$ then
$g({\bf z})$ is an eigenvalue of the matrix
$M_g:=[c_{j,s}]_{j,s=1}^s$ with corresponding eigenvector
$[b_1({\bf z}), \ldots,b_k({\bf z})]^T\neq 0$. Our assumption that
$\tilde{R}^{(k)}|_\varepsilon=0$ has rank $k$ implies that the
vectors $[b_1({\bf z}_s), \ldots,b_k({\bf z}_s)]^T $ for $s=1,
\ldots, k$ are linearly independent, thus they form a common
eigensystem for the matrices $M'_{x_1}|_{\varepsilon=0}, \ldots,
M'_{x_d}|_{\varepsilon=0}$. Thus, they pairwise commute and their
eigenvalues are the coordinates of ${\bf z}_1, \ldots, {\bf z}_k$,
proving the claim.

\end{proof}

\begin{rem}\emph{
Without further assumptions on the polynomials $b_1, \ldots, b_k$
we cannot guarantee that the polynomials
$f_{i,j}|_{\varepsilon=0}$ have no  roots outside of ${\bf z}_1,
\ldots, {\bf z}_k$. For example, if $k=d=1$ and ${\bf z}_1=c\neq
0$ but $b_1=x$, then $f_{11}=x^2-cx$ which also have $0$ as a
root. However, if we assume that $b_1, \ldots, b_k$ have no common
roots in $\mathbb{C}^m$ (e.g. $1\in \{b_1, \ldots, b_k\}$) then
all common roots of the polynomials $f_{i,j}|_{\varepsilon=0}$
correspond to eigenvalues and eigenvectors of
$M'_{x_i}|_{\varepsilon=0}$. Since ${\bf z_1}, \ldots, {\bf z_k}$
already provides a full system of eigenvectors for
$M'_{x_i}|_{\varepsilon=0}$, the polynomials
$f_{i,j}|_{\varepsilon=0}$ cannot have any other distinct root.}
\end{rem}

Our last result gives an asymptotic description of the roots of
the polynomials $\{f_{ij}\}$ in the case when $\varepsilon\neq 0$.
  Since the
coordinates of the vectors ${\bf v}_{i,j}$ are continuous in
$\varepsilon=0$  we can take their Taylor expansion around
$\varepsilon=0$ and consider them as elements  of the formal
series ring $\mathbb{C}[[\varepsilon]]$, as described in
Definition \ref{setting}. In this setting we will show that the
roots of the system  $\{f_{ij}\}$ are the centers of gravity (or
arithmetic means) of the clusters, modulo $\varepsilon^2$. Since
the arithmetic mean of a cluster is known to be better conditioned
than the individual roots in the clusters (c.f.
\cite{mandem95,corlessgiani97}), our result is therefore stable
for small enough values of $\varepsilon$.

\begin{prop}\label{meangen}
Let $B=[b_1,\ldots,b_n]$, $\{\bold{z}_1, \ldots, \bold{z}_k\}$ and
 for $i=1,\ldots, k$
\begin{scriptsize}
\[
\begin{aligned}
C_i=&\{[z_{i,1}+\delta_{i,1,1}\varepsilon,\ldots,z_{i,m}+\delta_{i,
1,m}\varepsilon], \ldots,\\
&\ldots[z_{i, 1}+\delta_{i,
n_i,1}\varepsilon,\ldots,z_{i,m}+\delta_{i,n_i,m}\varepsilon]\}\notag
\end{aligned}\]\end{scriptsize}
be as in Definition \ref{setting}. Let $\vec{\xi}_s=[\xi_{s,1},
\ldots, \xi_{s,m}]$ for $s=1, \ldots k$  be defined as
\begin{small}
\begin{equation}\label{xi}
 \xi_{s,i}:=
z_{s,i}+\frac{\sum_{r=1}^{n_s}\delta_{s,r,i}}{n_s}\varepsilon\quad
 i=1, \ldots, m.
\end{equation}
\end{small}
Then $\vec{\xi}_1, \ldots, \vec{\xi}_k$ satisfy modulo
$\varepsilon^2$ the defining equations $\{f_{i,j}\}$ of  the
approximate radical ideal of $C_1\cup \cdots \cup  C_k$ defined in
Definition \ref{apprad}.
\end{prop}

\begin{proof}
Fix $i\in \{1, \ldots, m\}$ and $j\in \{1, \ldots, k\}$. Define
$W$ to be the transpose of the $(k+1)\times k$ Vandermonde matrix
of $\vec{\xi}_1,\ldots, \vec{\xi}_k$ with respect to $[b_1,
\ldots, b_k,x_ib_{j}]$, i.e.
 $$
 W:=\left[b_t(\vec{\xi}_s)\left|(x_ib_{j})(\vec{\xi}_s)\right.\right]_{s,t=1}^k.
 $$ Also define $S$ to be the
$(k+1)\times k$ augmented matrix
 $$S:=\left[\left.\tilde{R}^{(k)}\right|{\bf r}_{i,j}\right]
 $$
 where $\tilde{R}^{(k)}$ and ${\bf r}_{i,j}$ was defined in Definition \ref{apprad}.
 Assume that
\begin{eqnarray}\label{inter}
x_ib_{j}-\left(\sum_{s=1}^k w_{i,j,s} b_s\right)\equiv 0 \mod
\varepsilon^2
 \end{eqnarray}
  is satisfied by $\vec{\xi}_s=[\xi_{s,1},
\ldots, \xi_{s,m}]$  for $s=1, \ldots, k$,  which is equivalent
for the column vector
 \begin{eqnarray}\label{cij} {\bf w}_{i,j}:=[w_{i,j,1}, \ldots,
w_{i,j,k}, -1]^T
\end{eqnarray}
to satisfy  the homogeneous linear system with coefficient matrix
$W$ modulo $\varepsilon^2$. On the other hand, from the definition
of the approximate radical ideal in Definition \ref{apprad}, we
also have that the augmented vector $[{\bf v}_{i,j}|-1]$  is a
solution of the homogeneous system corresponding to $S$. By our
assumption that $\det( \tilde{R}^{(k)}|_{\varepsilon=0})\neq 0$,
we also have that $\det (\tilde{R}^{(k)})\not \equiv 0\mod
\varepsilon^2$, which implies that both $S$ and $W$ have nullspace
of dimension $1$ modulo $\varepsilon^2$. Thus it is enough to show
that ${\bf w}_{i,j}$ is in the nullspace of $S$ modulo
$\varepsilon^2$, that will imply that ${\bf w}_{i,j}\equiv [{\bf
v}_{i,j}|-1]\mod \varepsilon^2$.

 Write
$$
{\bf w}_{i,j}\equiv {\bf w}_{i,j}^{(0)}+ {\bf
w}_{i,j}^{(1)}\varepsilon \quad W \equiv W^{(0)}+
W^{(1)}\varepsilon \quad S\equiv S^{(0)}+ S^{(1)}\varepsilon \mod
\varepsilon^2.
  $$
  At $\varepsilon =0$ we showed in the proof of Proposition \ref{continuous} that
  if
${\bf w}_{i,j}^{(0)}$ is in the nullspace of $W^{(0)}$ then it is
also in the nullspace of $S^{(0)}$.

It remains to prove that $ W^{(1)}{\bf w}_{i,j}^{(0)} +
W^{(0)}{\bf w}_{i,j}^{(1)}=0$ implies $ S^{(1)}{\bf w}_{i,j}^{(0)}
+ S^{(0)}{\bf w}_{i,j}^{(1)}=0$. We use the fact that
$$
S^{(0)}= V_1V_2^T \;\text{ and } \;S^{(1)}= \bar{V}_1 W^{(1)}+
\left(W_1^{(1)}\right)^T\bar{V}_2^T
$$
where $V_1$ and $V_2$ are the Vandermonde matrices of $C_1\cup
\cdots\cup C_k$ at $\varepsilon=0$ corresponding respectively  to
$[b_{\sigma(1)}, \ldots, b_{\sigma(k)}]$ and $[b_1, \ldots,
b_k,x_ib_{j}]$,   $W_1$ is the Vandermonde matrix corresponding to
$\vec{\xi}_1,\ldots, \vec{\xi}_k$ with respect to $(b_{\sigma(1)},
\ldots, b_{\sigma(k)})$, and $\bar{V}_1$ and $\bar{V}_2$ are the
same as $V_1$ and $V_2$, except the row corresponding to ${\bf
z}_s$ appears only once and it is multiplied by $n_s$. Since ${\bf
w}_{i,j}^{(0)}$ is in the nullspace of $W^{(0)}$, it is also in
the nullspace of $\bar{V}_2^T$, thus it remains to prove that
\begin{eqnarray}\label{simplify}
\bar{V}_1 W^{(1)}{\bf w}_{i,j}^{(0)}+ V_1V_2^T{\bf w}_{i,j}^{(1)}
=0.
\end{eqnarray}
Since  $W^{(1)}{\bf w}_{i,j}^{(0)}= - W^{(0)}{\bf w}_{i,j}^{(1)}$
by assumption,  (\ref{simplify}) is equivalent to
$$
\left[-\bar{V}_1 W^{(0)}+ V_1V_2^T\right] {\bf w}_{i,j}^{(1)} =0.
$$
But it is easy to see that $\bar{V}_1 W^{(0)}=V_1V_2^T$, which
proves the claim.

\end{proof}

As a corollary of the previous proposition we get that modulo
$\varepsilon^2$ the approximate multiplication matrices $M'_{x_1},
\ldots,
 M'_{x_d}$ form a pairwise commuting system of multiplication
 matrices for the roots $\vec{\xi}_1, \ldots, \vec{\xi}_k$.

\begin{cor}\label{commute}
Using the notation of Definition \ref{apprad} and Proposition
\ref{meangen} we have that for all $i=1, \ldots, k$ and $j=1,
\ldots, d$
 $$ M'_{x_j} {\bf e}_{\vec{\xi}_i} \equiv  \xi_{i,j}\; {\bf e}_{\vec{\xi}_i}
 \mod \varepsilon^2
 $$
 where
 $$
 {\bf e}_{\vec{\xi}_i}:=\left[b_s(\vec{\xi}_i)\right]_{s=1}^k.
 $$
 Thus the vectors $\{{\bf e}_{\vec{\xi}_i}\}_{i=1}^k$ form
 a common eigensystem  for the approximate multiplication matrices $M'_{x_1}, \ldots,
 M'_{x_d}$ modulo $\varepsilon^2$, which also implies that they
  are pairwise commuting
modulo $\varepsilon^2$, i.e. the entries of the commutators
$M'_{x_i}M'_{x_j}-M'_{x_j}M'_{x_i}$ are all divisible by
$\varepsilon^2$.
\end{cor}

\begin{rem}\emph{
In practice, for any particular choice of $\varepsilon \in
\mathbb{R}_+$ the system $\{ f_{i,j}\}$ is not necessary
consistent. Also, the approximate multiplication matrices
${M}'_{x_1}, \ldots, {M}'_{x_m}$ are not pairwise commuting, and
therefore not simultaneously diagonalizable.  However, one can
take any consistent subsystem of $\{ f_{i,j}\}$ such that it
defines each of the coordinates and solve this subsystem in order
to obtain the solutions. Another approach is the one described in
\cite{mandem95,corlessgiani97}:  If the distance of the clusters
from each other were order of magnitude larger than the size of
the clusters then a random linear combination of the matrices
$M'_{x_1}, \ldots, M'_{x_d}$ will have all its eigenvalues
distinct with high probablility. Using the eigensystem of this
random combination one can approximately diagonalize all of the
approximate multiplication matrices $M'_{x_1}, \ldots, M'_{x_d}$.
Then by Corollary \ref{commute} and \cite[Proposition
8]{corlessgiani97} the entries outside of the diagonal of the
resulting matrices will be small, asymptotically $\varepsilon^2$.
Taking the $i$-th diagonal entry of these nearly diagonal matrices
will give the coordinates of the $i$-th root of the approximate
radical, which by Proposition \ref{meangen} is approximately the
arithmetic mean of a cluster.}
\end{rem}

\begin{ex}\label{ex4}
Our last example is similar to Example \ref{ex3} but here we
increased the size of the clusters. Consider the polynomial system
given by {\small
\[
\begin{aligned}
\tilde f_1&=
x_1^2+3.99980x_1x_2-5.89970x_1+3.81765x_2^2-11.25296x_2\notag\\
&\qquad+8.33521\notag\\
\tilde f_2&=
x_1^3+12.68721x_1^2x_2-2.36353x_1^2+81.54846x_1x_2^2-177.31082x_1x_2\notag\\
&\qquad+73.43867x_1-x_2^3+6x_2^2+x_2+5\notag\\
\tilde f_3&=
x_1^3+8.04041x_1^2x_2-2.16167x_1^2+48.83937x_1x_2^2-106.72022x_1x_2\notag\\
&\qquad+44.00210x_1-x_2^3+4x_2^2+x_2+3\notag
\end{aligned}
\]}
which has a cluster of three common roots, $[0.8999, 1], [1, 1],
[1, 0.8999]$ around $[1,1]$ and a cluster of two common roots,
$[-1, 2], [-1.0999, 2]$ around $[-1,2]$. The clusters has size at
most $\varepsilon=0.1$. Using Chardin's subresultant method, we
obtained the multiplication matrices for this system, with respect
to the basis $B= [1,x_1,x_2,x_1x_2,x_1^2]$ and computed the matrix
of traces associated to the system, which is
\begin{scriptsize}
\begin{eqnarray}\label{rcmat}
R= \left[ \begin {array}{rrrrr} 5& 0.79999& 6.89990&
-1.40000& 5.01960\\
0.79999& 5.01960& -1.40000&
7.12928& 0.39812\\
6.89990& -1.40000& 10.80982&
-5.68988& 7.12928\\
-1.40000& 7.12928& -5.68988&
11.45876& -2.03262\\
5.01960& 0.39812& 7.12928& -2.03262& 5.11937
\end{array} \right]
\end{eqnarray}
\end{scriptsize}
After 2 steps of GECP on the matrix of traces we find the
partially reduced matrix $U_2$:
\begin{scriptsize}
\[
U_2=\left[ \begin {array}{rrrrr} 11.45876& -5.68988& 7.12928&
-1.40000& -2.03262\\
0& 7.98449& 2.14006& 6.20472& 6.11998\\
0& 0& 0.01039& 0.00799& 0.02243\\
0& 0& 0.00799& 0.00728& 0.01544\\
0& 0& 0.02243& 0.01544& 0.06796
\end{array} \right]
\]
\end{scriptsize}
with columns permuted to correspond to the basis $[x_1x_2, x_2,
x_1, 1, x_1^2]$.

 We also computed the singular values of
 $R$:
\[[24.06746,13.29215,
0.04397, 0.00362, 0.00035].\] We indeed have that the entries in
the last three rows of  $U_2$ and the third singular value
$\sigma_3$ are of the order of $\varepsilon^2$, which would
determine the numerical rank of $R$ to be 2.

By considering its last three rows of $U_2$  as zero, we compute
the nullspace of the resulting matrix, which gives the following
generators of $\widetilde{\sqrt{I}}/\tilde I$,
\begin{small}
\[
\begin{aligned} r_3&= x_2-1.46302+0.510803x_1  ,\notag\\
r_4&= x_1x_2+0.51920-1.505323x_1,\notag\\
r_5&= x_1^2-1.01587+0.08562x_1.\notag
\end{aligned}
\]
\end{small}
 From these we can define the multiplication matrices for $x_1$
and $x_2$ in $\mathbb{C}[\bold{x}]/\widetilde{\sqrt{I}}$ in the
basis $[1,x_1]$:
 {\scriptsize\[ \mathcal{M}'_{x_1}= \,
\left[ \begin {array}{cc}
 0& 1\\
1.01587& -0.08562
\end {array} \right] \quad
\mathcal{M}'_{x_2}= \, \left[ \begin {array}{cc}
1.46302& -0.51080\\
-0.51920& 1.50533\\
\end {array} \right]
\]}
These matrices do not commute but their commutator have small
entries:
 {\scriptsize$$  \mathcal{M}'_{x_1} \mathcal{M}'_{x_2}- \mathcal{M}'_{x_2} \mathcal{M}'_{x_1}=
    \left[ \begin {array}{cc}
 -0.000293 & -0.00143\\
 0.00147 & 0.000293
\end {array} \right] .
$$}
Thus the multiplication matrices are ``almost" simultaneously
diagonalizable. Following the method in \cite{corlessgiani97}, we
get the following approximate diagonalizations of
$\mathcal{M}'_{x_1}$ and $\mathcal{M}'_{x_2}$ using the eigenspace
of $\mathcal{M}'_{x_1} +\mathcal{M}'_{x_2}$:
 {\scriptsize$$
\mathcal{M}'_{x_1}\sim
 \left[ \begin {array}{cc}
                     -1.05162 & 0.001765\\
0.00116  &  0.966001\end {array} \right] \quad
\mathcal{M}'_{x_2}\sim \left[ \begin {array}{cc}
                  1.99959 & -0.001768\\
                  -0.001169  &   0.968759\end {array} \right].
$$}
The corresponding  diagonal entries give the solutions $[-1.05162,
1.99959]$  and $[0.966001, 0.968759]$ which are within $0.00167$
distance from the centers of gravity of the clusters in the
$\infty$-norm.
\end{ex}

\bibliographystyle{amsplain}
\providecommand{\bysame}{\leavevmode\hbox
to3em{\hrulefill}\thinspace}
\providecommand{\MR}{\relax\ifhmode\unskip\space\fi MR }
\providecommand{\MRhref}[2]{%
  \href{http://www.ams.org/mathscinet-getitem?mr=#1}{#2}
} \providecommand{\href}[2]{#2}

\end{document}